A Charging-as-a-Service Platform for Charging Electric Vehicles On the Move: New Vehicle Routing Model and Solution

**Jiahua Qiu**
Department of Civil and Coastal Engineering,
University of Florida, Gainesville, Florida, 32608
Email: jq22@ufl.edu

**Lili Du\***
Associate professor
Department of Civil and Coastal Engineering,
University of Florida, Gainesville, Florida, 32608
Email: lilidu@ufl.edu

**ABSTRACT**

Range anxiety has been the main challenge for the mass-adoption of electric vehicles (EVs). The emerging mobile electric-vehicle-to-electric-vehicle (mE2) charging technology, which allows an EV with the extra battery to charge the other EV on the move, offers a promising solution. Even though the physical feasibility of this new technology has been confirmed, we still face many operation problems. For example, how can we efficiently pair and route an electricity provider (EP) to a demand (ED) without introducing an extra detour? Motivated by this view, this study develops a Charging-as-a-Service (CaaS) platform, which seeks to dispatches the commercial EPs to serve the EDs, for cultivating this emerging charging service towards the low EV penetration environment. Mathematically, the CaaS platform is modeled as a vehicle routing problem (i.e., mE2-VRP), which optimally dispatches the EPs to approach and serve the EDs on the move while minimizing the EP fleet size and fulfilling all service requests. To adapt the CaaS platform to the online application in practice, we develop the Clustering-aided Clarke and Wright Savings (CCWS) algorithm to efficiently decompose and then solve the large-scale mE2-VRP by parallel computing. Our numerical experiments built upon citywide (Chicago) and statewide (Florida) found that the CCWS algorithm outperforms existing commercial solvers, and it enables us to investigate the performance of the CaaS platform under a realistic large-scale setup in a city or state. The CaaS performs better in low EV penetration markets while traffic congestion is mild and EDs require energy earlier. We can improve the performance by developing proper pricing strategies according to the EDs' energy requests and trip lengths.

# 1. Introduction

The transportation system accounted for the largest portion (28%) of greenhouse gas (GHG) emissions according to the data in the United States in 2018. The majority of which was related to the internal combustion engine vehicles (EPA, 2013). Governments around the world are setting ambitious targets for light carbon footprints. Battery-powered electric vehicles (EVs) are considered to be a promising solution since they will potentially reduce the oil dependence and the associated carbon emission. Accordingly, EVs are promoted by governments with different ranges of subsidies on both demand and supply sides (Foundation and Company, 2014). The global EVs sales exceeded 2 million, with 65% increase compared with that in 2017. However, existing studies (Gersdorf et al., 2020; Hao et al., 2020) showed that EV usage still stays at the early market stage. The main challenge that discourages the mass adoption of EVs is the 'range anxiety' resulted from the short driving range, long charging periods, and insufficient en route recharging services. The state-of-the-art shows that we are not able to enlarge battery capacity (Wu et al., 2020) in recent years. And it still takes about 30 minutes to fully charge the battery of an EV, even with supercharger stations, such as Level-3 charging stations or DC fast chargers (DCFC) (Roberts et al., 2017). This charging time is much longer than gas filling.

To relieve the range anxiety issue with current battery capacity, many solutions have been proposed, but each suffers from certain limits. The effectiveness has not been well confirmed. A brute force solution is to build sufficient charging stations to extend the driving range. But, this scheme cannot avoid the delay for charging. In addition, it always causes a high infrastructure investment for building charging stations and upgrading the existing power grid to avoid overload. Another solution is to build embedded wireless charging lanes. It enables charging on the move to save the charging time (He et al., 2018; Ngo et al., 2020), but is extremely expensive. Sweden has built a 1.24 electrified road, which cost $5.82 million (Karagiannopoulos, 2018). In view of the weakness of the above two solutions, some studies (Raeesi and Zografos, 2020; Shao et al., 2017) brought forward the idea of battery swapping (at a station or from other vehicles carrying extra batteries). However, EVs may use different models of batteries. This service hits the issues of battery ownership and incompatibility in practical implementation (Tang et al., 2020). The success of vehicle-to-vehicle charging technology in recent years raises other novel services (Cui et al., 2018b, 2018a; Huang et al., 2014; Li, 2019; Tang et al., 2020), which dispatch vehicles equipped with large batteries and chargers (i.e., mobile chargers (MC)) for delivering electricity to the low battery EVs at designated locations and times. Clearly, MCs may potentially reduce the 'range anxiety' to a certain extent, but travelers may still undergo delays resulting from the long charging time. Most recently, Connected and Autonomous Vehicle (CAV) initiates a breakthrough solution: mobile electric-vehicle-to-electric-vehicle (mE2) charging technique. Specifically, a pair of EVs (also CAVs) will cooperate their movement to "lock" their relative location for a period, so that they allow the on-the-move electricity exchange between them through either wireless power transfer (WPT) (Dai and Ludois, 2014; Hu et al., 2008) or converter-cable assembly charging technologies (Chakraborty et al., 2020; Roberts et al., 2017). As charging is conducted on the move, the mE2 charging technique is a promising solution to completely remove the range anxiety.

Even though the physical feasibility of this mE2 charging technology has been confirmed, the implementation of such charging service in the market still faces many operation issues. In particular, we still lack methods to quickly match the electricity suppliers (EPs) to the electricity demands (EDs), and then efficiently route the EPs to approach the EDs, considering their respective mobility and trip plans. Only a few studies investigated these issues in the literature. The studies of (Kosmanos et al., 2018; Maglaras et al., 2014) developed a mobile charging system, in which the EPs are buses that carry high-capacity batteries. The EDs need to move with a bus while receiving energy via plug-in electric connection or via WPT (electromagnetic induction). Clearly, the fixed bus routes can only provide limited service coverage.

Another pioneering study developed by (Abdolmaleki et al., 2019) tried a crowdsourcing approach, by which a mathematical program was developed to route the peer volunteer suppliers with the objective to maximize the number of served electricity demands. The experiment results show that each electricity demand vehicle needs more than three volunteer suppliers on average to complete a trip. The crowdsourcing approach may meet the difficulty to find enough volunteer suppliers in the early market with a limited number of EVs.

To address the above issues, this study seeks to develop a Charging-as-a-Service (CaaS) platform, which takes advantage of the mE2 technology and is able to cultivate this service when the EV market penetration is still low. It is a centralized platform, which optimally routes the commercial EPs (also CAVs) fleet to deliver electricity on the move for the EDs (also CAVs) so that the EDs don't need to detour or wait for the charging. More exactly, the platform collects the EDs' charging requests and trip plans within a lead time, and then finds optimal routes to dispatch the EPs for the EDs with the objective to provide timely charging services en route. Considering the demand comes ad hoc, the online service decision process can be conducted batch by batch. In addition, it is aware that these EPs introduce extra vehicles on the roads and may cause severe traffic congestion. To mitigate this traffic overhead, we want to limit the fleet size of the EPs, which in turn will affect the service efficiency.

With the above considerations and dilemma, we mathematically model the service of the CaaS platform in each batch as a vehicle routing problem (VRP) for mobile EV-to-EV charging (labeled as mE2-VRP). It aims to efficiently route the EPs to approach and then deliver electricity to the EDs timely, while they are moving along the pre-defined paths. The mE2-VRP is different from existing VRP models in several aspects. First of all, the mE2-VRP considers mobile demands and on-the-move delivery while most of the existing VPR models consider static delivery/demand objects. From this view, the most relevant model is the VRP with roaming delivery locations (VRPRDL) investigated by (Ozbaygin et al., 2017; Reyes et al., 2017). Specifically, the VRPRDL considers that each customer will visit a set of static locations, e.g., home and workplaces, and a single delivery can be at one of the static candidate locations. However, the mE2-VRP will have an EP-ED pair run together so that the electricity delivery can be conducted on the move. This requires the mE2-VRP to synchronize the movement of an EP-ED pair in a spatiotemporal space rather than a snapshot at a static location. Moreover, the delivery product of the mE2-VRP is electricity, which can be requested by an ED at any time with various amounts during its route. This unique feature makes an ED can accept multiple times of services from different EPs during its route (i.e., partial charging is allowed at each service). Last, the CaaS platform raises new computation and scalability difficulties. Specifically, the CaaS platform will serve a large number of EDs running en routes in a local network. Therefore, we need to solve the mE2-VRP quickly so that we can provide timely services before EDs run out of battery en routes. On the other hand, it is often very challenging to solve a large-scale VRP. All the aforementioned unique features make the mE2-VRP call for new formulations and solution methodology, which thus motivate the research efforts in this paper.

The main efforts and methodology contributions of this study are summarized as follows. First of all, this study introduces a commercial CaaS platform leveraging the emerging EV-to-EV charging and CAV technologies. This success of this CaaS platform is able to relieve EV range anxiety and cultivate the mass adoption of EVs without prohibitive infrastructure investment. Second, this CaaS is modeled as a novel mE2-VRP model, which optimally routes EPs to charge EDs *on the move* while minimizing the traffic overhead introduced by the EP fleet and maximizing the service efficiency of the CaaS. Two critical points highlight the novelty here. The mE2-VRP model is built upon an encounter network and addresses trip synchronization difficulty among multiple mobile vehicles. The mE2-VPR develops new constraints to eliminate the sub-tours uniquely introduced by the on-the-move electricity delivery. Their effectiveness is mathematically proved. Moreover, we develop a Clustering aided Clarke and Wright Savings (CCWS) algorithm to address the scalability difficulties so that we can implement the CaaS in practice. Mainly, we first contribute a quantitative approach to measure the potential that two EDs can be served by the same EP along their trips. Based on this measurement, we develop a customized clustering method that decomposes a large-scale mE2-VRP into multiple small subproblems. And then, we solve the subproblems efficiently and independently by parallel computing to obtain a feasible solution of the mE2-VRP. Next, we design a merging algorithm that improves the feasible solution to a better solution for the mE2-VRP. The CCWS adapts the CaaS platform to the online application involving large-scale ad hoc demands in practice.

Taking advantage of the mE2-VRP model and the CCWS algorithm, we conducted numerical studies to explore the insights and the applicability of such platform in practical implementations. For example, by analyzing the traffic overhead introduced by the EPs fleet, we investigated the adoption of this CaaS platform under different EV market penetrations and congestion levels. The results imply that the CaaS platform is applicable for the next few decades according to the projected low EV market penetrations. In addition, proper pricing strategies according to the EDs' energy requests and trip lengths will help improve the service efficiency and ultimately relieve travelers' range anxieties for using EVs for their trips. Note that without an efficient solution approach like the CCWS algorithm, we are not able to investigate those insights for the CaaS implemented in a city-wide or state-wide network.

Overall, this study is among the first efforts to study the application of the on-the-move vehicle-to-vehicle charging technologies for EV charging services in practice. It significantly contributes to both the methodology development in literature and applications in practices for promoting electrification in transportation. In addition, the mE2-VRP model and the CCWS algorithm potentially enable more advanced collaborative mobile delivery services in the future. For example, by using drones, a package can be optimally relayed by several trucks on the move along their trips until arriving its destination. Here, the mE2-VRP can optimally schedule the trips of drones so that they can deliver packages from one truck to another on the move. By doing so, we can fully use the idle capacities of trucks on the road and also save their trips.

The efforts of this study are organized by the following structure in the rest of the paper. Section 2 introduces the CaaS platform. Following that, we formulate the mE2-VRP in Section 3 and develop the solution methodology in Section 4. Based on the Chicago sketch network and Florida statewide network, Section 5 conducts numerical experiments to validate the developed models and algorithms. Some interesting managerial insights are also presented. Last, we summarize the whole study and propose the potential future work in Section 6.

## 2. Problem description and formulation

The CaaS platform seeks to efficiently provide an on-the-move charging service for the electric vehicles en routes. Before formally formulate this problem, we would first clarify several assumptions. First of all, the CaaS platform does not require an ED to detour for accepting the electricity delivery. Instead, an ED may meet an EP and accept the electricity delivery on the move along the ED's pre-defined route. It may cause a waiting time before an ED meets the assigned EP at a location on its route. The CaaS will satisfy each ED's tolerable delay reported in a service request. Therefore, the CaaS will help an ED avoid all other unwanted delays due to a recharging service. Moreover, by taking advantage of the mE2 charging technique, the CaaS platform allows an ED to accept multiple times of services from different EPs at different spots along its route. This means that partial charging at each electricity delivery service is permitted. Based on these assumptions, we will explore optimal routes for dispatching the EPs so that they can serve the EDs en route with a minimum fleet size, which helps reduce the traffic overhead introduced by the EPs and also brings benefit to save the operational cost[1]. To do that, the sections below will formally define the CaaS platform, electricity suppliers, and demands.

## 2.1 CaaS platform and electricity suppliers

To model the CaaS platform, we first formally describe a mE2 service, taking Figure 1(a) as an example. Specifically, we consider the platform can serve ED(i) at the earliest time $t_0^i$ at the location $i_0$ on its route. After passing location $i_0$, ED(i) goes through node $n_1$ at the time $t_{n_1}^i$ and then arrives at destination $i_e$ at the time $t_e^i$. In the meantime, an EP departs from its current location or depot and will serve an ED such as $i$ by meeting it at a node such as $i_0$ and then move together with ED(i) to distribute electricity along its route. After finishing with ED(i), the EP will continue to serve another ED if it still has enough electricity inventory. By studying the chance that an EP moves from one ED to another, we notice that those nodes where two EDs encounter (such as $n_1$), are more attractive to schedule the mE2 service since it potentially enables an EP to serve a new ED (such as $j$) without an extra idle travel after it finishes with the current ED (such as $i$). This unique feature encourages this study to build the CaaS platform on a novel network mainly composed of the encounter nodes of EDs' routes. Then, an EP trip formed in such a network will enable this EP to meet and serve several EDs at the encounter nodes along its route. Thus, this novel network will help us to synchronize the trips of EPs and EDs.

Consequently, we consider a directed graph $G(\mathcal{N}, A)$ with the node set $\mathcal{N} = N \cup N_0 \cup N_e \cup N_p$ and the arc set $A$ (see Figure 1(b) for an example). Specifically, $N$ includes all the encounter nodes between the EDs. $N_0$ represents the set of the first available locations along the itineraries of the EDs, where the platform can provide the earliest mE2 charging service. $N_e$ represents the set of the destinations of the EDs. $N_p$ contains the current locations of the available EPs (can be the depots or the current locations of the idle EPs). Each arc $(n,m) \in A, \forall n, m \in \mathcal{N}, n \neq m$ in the graph $G(\mathcal{N}, A)$ represents a route connecting two nodes in the physical map. Note that a physical route may involve multiple roads which are not presented in the graph $G(\mathcal{N}, A)$. Each arc in $G$ has an associated average travel time $t_{n \to m}$ and energy cost $e_{n \to m}$. The visiting time stamps of each ED at the nodes on its route are recorded by a set of continuous variables and we will introduce them later.

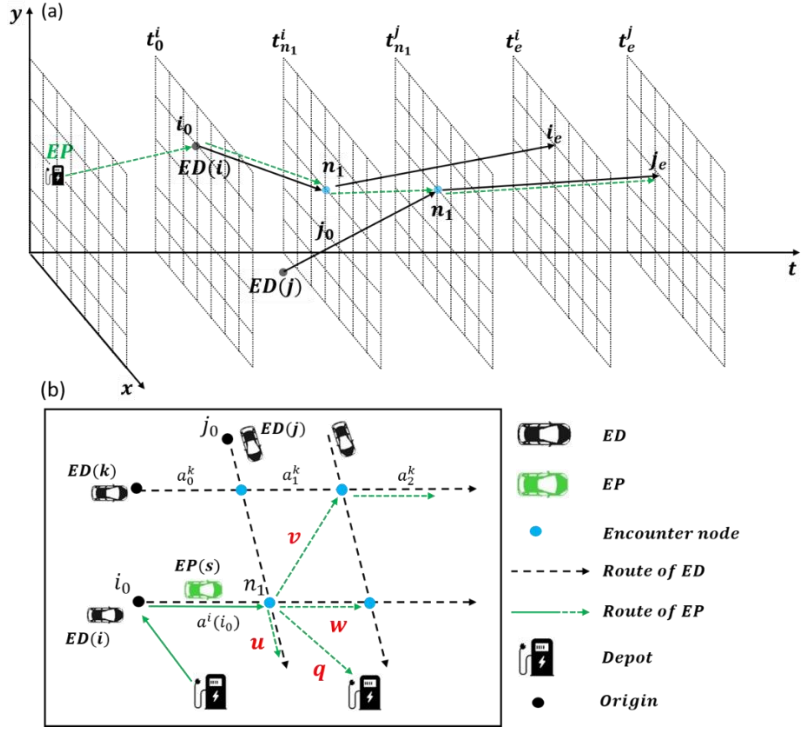

**Figure 1.** Schematic representation of the mE2-VRP on a (a) spatiotemporal network and (b) simplified spatial network

---

[1] We consider average energy consumption and travel time on each link, which is a well-accepted assumption in the literature (Bongiovanni et al., 2019; Cui et al., 2018b, 2018a; Raeesi and Zografos, 2020; Schneider et al., 2014; Tang et al., 2020). The potential mismatch between the EPs and the EDs in the routing plan resulting from this assumption can be accommodated by the demands' schedule flexibility (discussed in Section 2.2).

Next, we consider a set of EPs denoted by $\mathcal{S}$. They are either parking in the depots or idle en routes after completing their services in the last batch. For each EP($s$), $\forall s \in \mathcal{S}$, we denote its initial energy inventory by $e_0^s$, its current location by $s_0$, and its energy consumption rate by $\varpi_s$. We also let $e_{n \to m}^-, \forall n, m \in \mathcal{N}$ represent the average energy consumption for an EP traveling on the path from node $n$ to $m$ involving several arcs in $G(\mathcal{N}, A)$, which can be estimated by the EP's energy consumption rate and the average route travel time. This study discretizes the service horizon, i.e., one day, into $H$ time intervals, with interval length $T$. Accordingly, we denote $\mathcal{D}_h, \forall h \in H$, as a set of the EDs that send the mE2 charging service requests within the $h$-th time interval. We further assume that at the beginning of each time interval, $h$, the total energy storage of the available EPs is sufficient to serve all energy requests. The platform solves the mE2-VRP once in each time interval and then dispatch a fleet of the EPs to serve the EDs in $\mathcal{D}_{h-1}$. For convenience, we abbreviate $\mathcal{D}_h$ as $\mathcal{D}$ for the rest of the paper since the approach is adaptive to each time interval $h \in H$. As the platform is operated under a rolling horizon, the solution time of the mE2-VRP is required to be smaller than the interval length $T$, which calls for an efficient solution algorithm.

## 2.2 Electricity demands en route

We next define the features of the electricity demands. When the CaaS platform receives the electricity requests from the EDs, it takes a tolerant service lead time to process the requests and dispatch the EP services. Accordingly, ahead of the leading time, each ED(i) is required to provide 1) its battery features (capacity $\bar{e}_i$, and energy consumption rate $\varpi_i$), and 2) its origin-destination ($i_0, i_e$) and route information. Here, the origin $i_0$ is the first location[2] that the ED is available for accepting the service. Besides, the system needs information regarding its estimated arrival time ($\tilde{t}_0^i$), battery level ($e_0^i$), and endurable waiting time $\bar{\tau}_i$ at its origin $i_0$.

With the given information of the EDs, the system re-defines the route of each ED on the graph $G(\mathcal{N}, A)$ by a tuple of nodes, i.e., $\mathcal{N}_i = (i_0, N_i, i_e), i \in \mathcal{D}$, where $N_i$ represents the encounter nodes along the route of ED (i). Moreover, we use $t_n^i$ to represent the departure time of ED(i) at node $n \in \mathcal{N}_i$. Assuming each ED(i) is only willing to wait for the service at the origin $i_0$, we can estimate its departure time from the location $i_0$ by $t_0^i = \tilde{t}_0^i + \tau_i$, where the waiting time $\tau_i$ satisfies $0 \leq \tau_i \leq \bar{\tau}_i$. Accordingly, along the route of ED(i) on the graph $G(\mathcal{N}, A)$, we can have the relation $t_n^i = t_0^i + t_{i_0 \to n} = \tilde{t}_0^i + t_{i_0 \to n} + \tau_i, \forall i \in \mathcal{D}, n \in \mathcal{N}_i$, where $t_{i_0 \to n}$ represents the average travel time from location $i_0$ to node $n$ along the route on $G(\mathcal{N}, A)$. A set of the EDs passing through a node $n$, is denoted by $D_n, n \in \mathcal{N}$.

To facilitate the development of our mathematical model, this study also presents the route of each ED(i) on $G(\mathcal{N}, A)$ by the arcs involved. Accordingly, we let $A_i = \{a_n^i\}_{n \in \mathcal{N}_i} \subseteq A$ denote the route of ED(i), where $a_n^i$ represents an arc starting from node $n$ on the route of ED(i). To simplify the presentation, $a_n^i$ is replaced by $a^i(n)$ in some formulations. Figure 1(b) demonstrates the examples, in which two different ways are used to denote the arc for ED(i) and ED(k) respectively. We further introduce continuous variables $e_{i,a}^+, i \in \mathcal{D}, a \in A_i$ to represent the electricity delivered to ED(i) on the arc $a$. Correspondingly, we use parameters $e_{i,a}^-, i \in \mathcal{D}, a \in A_i$ to represent the energy loss of ED(i) along the arc $a$. Note that we consider $e_{i,a}^-$ as a known parameter that can be estimated by the given terrain type, distance, and average travel time.

---

[2] The first service location, $i_0$ can any spot en route or even at the origin if the ED does not depart.

## 2.3 Encountering and Switching Schemes

This study considers that an EP(s) can meet an ED(j) at either its origin node or any encounter nodes, and then follows ED(i) on the arc $a \in A_i$ to deliver the electricity until the next node, with known the power transfer rate $\eta$. If EP(s) has enough energy once it arrives at the next node, EP(s) can continue to charge ED(i); or it can either switches to serve another ED(j) at the current node (local customer switch, shorten it as local switch hereafter) or drives to serve another ED(k) at different node (distant customer switch, shorten as distant switch hereafter). If EP(s) does not have enough energy, it will return to a depot. To model these decision actions involved in this mE2 charging service, we introduce a set of binary planning variables as follows.

$$z_{i,a}^s = \begin{cases} 1 & \text{if } EP(s) \text{ serves } ED(i) \text{ on the arc } a \in A_i \\ 0 & \text{otherwise} \end{cases}, \quad \forall s \in \mathcal{S}, i \in \mathcal{D}, a \in A_i$$

$$o_i^{n,s} = \begin{cases} 1 & \text{if } EP(s) \text{ travels from a depot} \\ & \text{to charge the } ED(i) \text{ at node } n \\ 0 & \text{otherwise} \end{cases}, \quad \forall s \in \mathcal{S}, i \in \mathcal{D}, n \in \mathcal{N}_i \backslash i_e$$

$$q_i^{n,s} = \begin{cases} 1 & \text{if } EP(s) \text{ leaves } ED(i) \text{ and} \\ & \text{return to depot from node } n, \\ 0 & \text{otherwise} \end{cases}, \quad \forall s \in \mathcal{S}, i \in \mathcal{D}, n \in \mathcal{N}_i \backslash i_0$$

$$u_{i,j}^{n,s} = \begin{cases} 1 & \text{if } EP(s) \text{ conducts local switches} \\ & \text{from } ED(i) \text{ to } ED(j) \text{ at node } n \\ 0 & \text{otherwise} \end{cases}, \quad \forall s \in \mathcal{S}, i \in \mathcal{D}, n \in \mathcal{N}_i \backslash i_0, j \in D_n, j \neq i$$

$$v_{i,j}^{n,m,s} = \begin{cases} 1 & \text{if } EP(s) \text{ conducts distant switches from} \\ & ED(i) \text{ at node } n \text{ to } ED(j) \text{ at node } m \\ 0 & \text{otherwise} \end{cases}, \quad \forall s \in \mathcal{S}, i,j \in \mathcal{D}, i \neq j, n \in \mathcal{N}_i \backslash i_0, m \in \mathcal{N}_j \backslash j_e, n \neq m$$

$$w_i^{n,s} = \begin{cases} 1 & \text{if } EP(s) \text{ charges } ED(i) \text{ on arc } a^i(n) - 1 \\ & \text{and continue to charge } ED(i) \text{ on arc } a^i(n) \\ 0 & \text{otherwise} \end{cases}, \quad \forall s \in \mathcal{S}, i \in \mathcal{D}, n \in \mathcal{N}_i$$

Based on the above problem setup, we try to find optimal schemes to match and route the EPs to the EDs so that they meet and conduct the electricity charging when they are on the move. For example, in Figure 1(a), EP(s) firstly charges ED(i) from node $i_0$ to node $n_1$. Then EP(s) waits $t_{n_1}^j - t_{n_1}^i$ units of time at node $n_1$, switches the charging target to ED(j) and continues to charge ED(j) until its destination $i_e$. To generate an optimal schedule for the CaaS platform to dispatching the EPs, we developed a mE2-VRP model, which is a large-scale mixed-integer linear programming (MILP), and solved it by a heuristic solution approach named as CCWS algorithm to adapt this online application. We discuss the technical details of these methods in the following sections. To facilitate the understanding of our methods, we summarize the variables and parameters shown in **Table 2** in Appendix A.

## 3. Mathematical formulation

The service of the CaaS platform raises a VRP problem for the electricity delivery between two mobile electrical vehicles (i.e., mE2-VRP). Mathematically we formulated the mE2-VRP as a MILP from Equation (1) to Equation (15), which explores the optimal routing plans for the EPs to serve a set of EDs subject to their energy safety and charging requests, while minimizing the platform operation cost and extra traffic introduced by EPs (i.e., traffic overhead). Accordingly, the objective function in Equation (1) seeks to minimize the fleet size of the EPs. The development of the constraints considers the following aspects: flow conservation of the EPs ensured by the constraint set (2)- (7); the spatial- and energy- feasibility of each local and distant switch ensured by the constraints (8)- (12). Last, constraints (13) and (14) enable the temporal-feasibility of a local/distant switch and they eliminate the illegal subtours of the EPs. We introduce the thoughts to develop those constraints in detail by the following subsections.

**mE2-VRP**

$$\min F(O) = \sum_{s \in \mathcal{S}} \sum_{i \in \mathcal{D}} \sum_{n \in \mathcal{N}_i \setminus i_e} o_i^{n,s} \tag{1}$$

subject to

$$\sum_{\substack{j \in D_n \\ j \neq i}} u_{i,j}^{n,s} + \sum_{\substack{j \in \mathcal{D} \\ j \neq i}} \sum_{\substack{m \in \mathcal{N}_j \setminus j_e \\ n \neq m}} v_{i,j}^{n,m,s} + w_i^{n,s} + q_i^{n,s} = z_{i,a^i(n)-1}^s, \qquad \forall s \in \mathcal{S}, i \in \mathcal{D}, n \in N_i \tag{2}$$

$$\sum_{\substack{j \in D_n \\ j \neq i}} u_{j,i}^{n,s} + \sum_{\substack{j \in \mathcal{D} \\ j \neq i}} \sum_{\substack{m \in \mathcal{N}_j \setminus j_0 \\ m \neq n}} v_{j,i}^{m,n,s} + o_i^{n,s} + w_i^{n,s} = z_{i,a^i(n)}^s, \qquad \forall s \in \mathcal{S}, i \in \mathcal{D}, n \in N_i \tag{3}$$

$$\sum_{\substack{j \in D_{i_e} \\ j \neq i}} u_{i,j}^{i_e,s} + \sum_{\substack{j \in \mathcal{D} \\ j \neq i}} \sum_{\substack{m \in \mathcal{N}_j \setminus j_e \\ m \neq i_e}} v_{i,j}^{i_e,m,s} + q_i^{i_e,s} = z_{i,a^i(i_e)-1}^s, \qquad \forall s \in \mathcal{S}, i \in \mathcal{D} \tag{4}$$

$$\sum_{\substack{j \in D_{i_0} \\ j \neq i}} u_{j,i}^{i_0,s} + \sum_{\substack{j \in \mathcal{D} \\ j \neq i}} \sum_{\substack{m \in \mathcal{N}_j \setminus j_0 \\ m \neq i_0}} v_{j,i}^{m,i_0,s} + o_i^{i_0,s} = z_{i,a^i(i_0)}^s, \qquad \forall s \in \mathcal{S}, i \in \mathcal{D} \tag{5}$$

$$\sum_{i \in \mathcal{D}} \sum_{n \in \mathcal{N}_i \setminus i_e} o_i^{n,s} \leq 1, \qquad \forall s \in \mathcal{S} \tag{6}$$

$$\sum_{i \in \mathcal{D}} \sum_{n \in \mathcal{N}_i \setminus i_0} q_i^{n,s} \leq 1, \qquad \forall s \in \mathcal{S} \tag{7}$$

$$e_n^i = e_0^i + \sum_{a=a^i(i_0)}^{a^i(n)-1} e_{i,a}^+ \left( \sum_{s \in \mathcal{S}} z_{i,a}^s \right) - \sum_{a=a^i(i_0)}^{a^i(n)-1} e_{i,a}^-, \qquad \forall i \in \mathcal{D}, n \in \mathcal{N}_i \setminus i_0 \tag{8}$$

$$0 \leq e_{i,a}^+ \leq \bar{e}_{i,a}^+, \qquad \forall i \in \mathcal{D}, a \in A_i \tag{9}$$

$$\underline{e}_d \leq e_n^i \leq \bar{e}_i, \qquad \forall i \in \mathcal{D}, n \in \mathcal{N}_i \setminus i_0 \tag{10}$$

$$E_s = \sum_{i \in \mathcal{D}} \sum_{n \in \mathcal{N}_i \setminus i_e} o_i^{n,s} e_{s_0 \to n}^- + \sum_{i \in \mathcal{D}} \sum_{\substack{j \in \mathcal{D} \\ j \neq i}} \sum_{\substack{n \in \mathcal{N}_i \setminus i_0 \\ m \in \mathcal{N}_j \setminus j_e}} v_{i,j}^{n,m,s} e_{n \to m}^- + \tag{11}$$

$$\sum_{i \in \mathcal{D}} \sum_{a=a^i(i_0)}^{a^i(i_e)-1} z_{i,a}^s (e_{i,a}^+ + e_{s,a}^-) + \sum_{i \in \mathcal{D}} \sum_{n \in \mathcal{N}_i \setminus i_0} q_i^{n,s} e_{n \to p_n}^-$$

$$e_0^s - E_s \geq \underline{e}_s, \qquad \forall s \in \mathcal{S} \tag{12}$$

$$t_n^j - t_n^i \geq M(u_{i,j}^{n,s} - 1), \qquad \forall s \in \mathcal{S}, i \in \mathcal{D}, n \in \mathcal{N}_i \setminus i_0, j \in \mathcal{D}_n, j \neq i \quad (13)$$

$$t_m^j - t_n^i - t_{n \to m} \geq M(v_{i,j}^{n,m,s} - 1), \qquad \forall s \in \mathcal{S}, i,j \in \mathcal{D}, i \neq j, c \in \mathcal{N}_i \setminus i_0, m \in \mathcal{N}_j \setminus j_e \quad (14)$$

$$t_n^i = t_0^i + t_{i_0 \to n} = \tilde{t}_0^i + t_{i_0 \to n} + \tau_i, \qquad \forall i \in \mathcal{D}, n \in \mathcal{N}_i \quad (15)$$

### 3.1 Flow conservation constraints

We first introduce the thoughts of developing the constraint set (2)-(5) to ensure the EP flow conservation at each node. Mainly, constraints (2) and (4) model the actions that an EP will take at a node $n$. To do that, we introduce binary variables $z = \{z_{i,a}^s\}$ to indicate if there is an EP moving on the arc $a$. Taking Figure 2(a) as an example, we can see that EP($s$) charges ED($i$) as they move together on arc $a^i(n_1) - 1$ and arrives at node $n_1$. Thus, there is an EP entering node $n_1$. Mathematically, it indicates that $z_{i,a^i(n_1)-1}^s = 1$ on the right side of (2). To ensure the flow conservation, we should have an EP leaving node $n_1$. Figure 2(a) shows that there are four options, and EP($s$) can only take one out of them. This study thus uses four binary variables on the left-hand side of the constraint (2) to model this decision. Specifically, EP($s$) may continue charging ED($i$) toward the downstream segment, then we have $w_i^{n_1,s} = 1$; or do a local service switch from ED($i$) to ED($k$) at node $n_1$, then we have $u_{i,k}^{n_1,s} = 1$; or do a distant service switch from ED($i$) to ED(j) by traveling from node $n_1$ to node $n_3$, then we have $v_{i,j}^{n_1,n_3,s} = 1$; or return to a depot, which means $q_i^{n_1,s} = 1$. However, if $z_{i,a^i(n_1)-1}^s = 0$ on the right side of (2), none of the four actions are feasible for EP($s$) at node $n_1$. The flow conversation constraint (2) makes the four integer variables on the left side take zero values. Similarly, constraint (4) models the same flow conservation relationship in the case where node $n_1$ is the destination of ED($i$), $i_e$. Note that variable $w_i^{n_1,s}$ is omitted in constraint (4) because ED($i$) arrives at destination and there is not a downstream arc on its route.

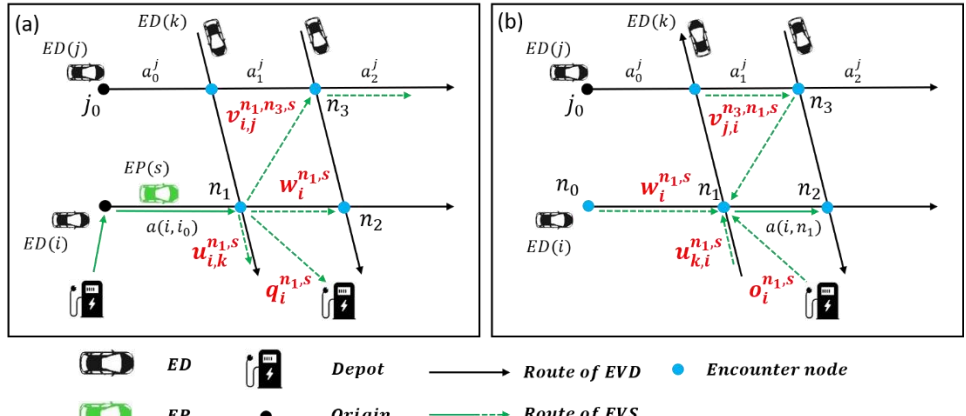

**Figure 2.** Representation of the flow conservation constraints.

Similarly, constraints (3) and (5) capture all possible actions that drive an EP arriving at node $n$. Figure 2(b) shows an example. If EP($s$) serves ED($i$) on the arc $a^i(n_1)$, then we know that there must be an EP leaving from node $n_1$, which makes $z_{i,a^i(n_1)}^s = 1$ on the right side of (3). On the other hand, the left side of (3) indicates that EP($s$) can arrive at node $n_1$ by only taking one of the four actions, which introduces four integer variables similar to what we discussed above. Namely, it may come from a depot, then we have an integer variable $o_i^{c_1,s} = 1$, or continue charging ED(i) from an upstream segment to node $n_1$, then it makes $w_i^{n_1,s} = 1$; or perform a local switch from ED($k$) to ED($i$) at $n_1$, then we have $u_{k,i}^{n_1,s} = 1$; or conduct a distant switch from ED(j) to ED(i), then we have $v_{j,i}^{n_3,n_1,s} = 1$. Clearly, if $z_{i,a^i(n_1)}^s = 0$ on the right side of (3), the constraint makes the corresponding four integer variables on the right side take zero values since there is not an EP arriving at node $n_1$. Similarly, this flow conservation relationship is modeled by constraint (5) when $n_1$ is the origin of ED($i$), $i_0$. In this case, there is no upstream arc to $i_0$ on the route of ED($i$), so variable $w_i^{n_1,s}$ is omitted in constraint (5).

Furthermore, constraint (6) ensures that an EP will only be assigned to one ED at most when it departs from a depot. Constraint (7) ensures that an EP can return to a depot at most once within a batch of service.

## 3.2 ED and EP energy constraint

We next discuss constraints (8)- (12) to demonstrate how our model ensures the feasible energy inventory of the EDs and the EPs during the service. First of all, the constraint set (8)- (10) ensures the batteries of the EDs are not depleted during the trips. Specifically, $e_n^i$ is used to denote the energy inventory of ED($i$) at any node $n \in \mathcal{N}_i \setminus i_0$. It is formulated by three items in constraint (8). The first item of constraint (8) is the energy inventory of ED($i$) upon departure. The second item of constraint (8) is the total energy that ED(i) receives before arriving at node $n$, where $e_{i,a}^+ (\sum_{s \in \mathcal{S}} z_{i,a}^s)$ indicates the electricity that ED(i) receives along the arc $a$. The last item of Equation (8) is the total energy that ED(i) consumes before it arrives at node $n$. Note that the energy reception in the second item of Equation (8) is non-linear. Therefore, to facilitate computation, we introduce auxiliary variables $\gamma_{i,a} = e_{i,a}^+ (\sum_{s \in \mathcal{S}} z_{i,a}^s), \forall i \in \mathcal{D}, a \in A_i$ and develop Equations (8.1)- (8.4) below to linearize the non-linear term of $e_{i,a}^+ (\sum_{s \in \mathcal{S}} z_{i,a}^s)$ in constraint (8), where $M$ is a sufficient big value.

$$\gamma_{i,a} - e_{i,a}^+ \leq 0, \ \forall i \in \mathcal{D}, a \in A_i \tag{8.1}$$

$$\gamma_{i,a} - M\left(\sum_{s \in \mathcal{S}} z_{i,a}^s\right) \leq 0, \ \forall i \in \mathcal{D}, a \in A_i \tag{8.2}$$

$$M\left(\sum_{s \in \mathcal{S}} z_{i,a}^s - 1\right) - \gamma_{i,a} + e_{i,a}^+ \leq 0, \ \forall i \in \mathcal{D}, a \in A_i \tag{8.3}$$

$$\gamma_{i,a} \geq 0, \ \forall i \in \mathcal{D}, a \in A_i \tag{8.4}$$

More exactly, Equations (8.1) and (8.3) ensure $\gamma_{i,a} = e_{i,a}^+$ when $\sum_{s \in \mathcal{S}} z_{i,a}^s = 1$. Equations (8.2) and (8.4) ensures $\gamma_{i,a} = 0$ when $\sum_{s \in \mathcal{S}} z_{i,a}^s = 0$. Constraint (9) bounds the variable $e_{i,a}^+$ with given maximum energy that ED(i) can receive on arc $a$ (i.e., parameter $\bar{e}_{i,a}^+$). It depends on the travel time to go through the arc and power transfer rate $\eta$. Constraint (10) ensures that the battery inventory of the EDs is always above safety inventory $\underline{e}_d$ and not exceeds battery capacity $\bar{e}_i$.

On the other hand, the routing plan will ensure the energy safety so that each EP can return to a depot for recharging. We model the total energy loss of an EP during a trip, $E_s$, by Equation (11). It includes four items respectively counting the energy loss of an EP before, during and after the charging service in a trip. More exactly, the first item ($\sum \sum o_i^{n,s} e_{s_0 \to n}^-$) represents the idle energy loss for an EP to reach the first service reception location $n$ from its initial departure location $s_0$ (i.e., $o_i^{n,s} = 1$). The second item ($\sum \sum \sum v_{i,j}^{n,m,s} e_{n \to m}^-$) is the required energy for an EP to conduct distant switch services over all possibilities. The third item ($\sum \sum z_{i,a}^s e_{i,a}^+ + \sum \sum z_{i,a}^s e_{s,a}^-$) consists of the energy distributed to EDs ($z_{i,a}^s e_{i,a}^+$) and self-consumed ($z_{i,a}^s e_{s,a}^-$) during the services, and the last item ($\sum \sum q_i^{n,s} e_{n \to p_n}^-$) represents the idle energy loss for an EP returning to the nearest depot $p_n$ from node $n$. Built upon constraint (11), constraint (12) ensures that the energy inventory of an EP is always above the safety inventory $\underline{e}_s$ during the service, where $e_0^s$ is the initial energy inventory of an EP. Moreover, we use the similar linearization approach adopted in the constraint set (8.1)- (8.4) to linearize the non-linear term $z_{i,a}^s e_{i,a}^+$ in Equation (11). With additional set of auxiliary variables, $\zeta_{i,a}^s = z_{i,a}^s e_{i,a}^+, \forall i \in \mathcal{D}, a \in A_i, s \in \mathcal{S}$, constraints set (11.1) to (11.4) are provided below.

$$\zeta_{i,a}^s - e_{i,a}^+ \leq 0, \quad \forall i \in \mathcal{D}, a \in A_i, s \in \mathcal{S} \tag{11.1}$$

$$\zeta_{i,a}^s - M z_{i,a}^s \leq 0, \quad \forall i \in \mathcal{D}, a \in A_i, s \in \mathcal{S} \tag{11.2}$$

$$M(z_{i,a}^s - 1) - \zeta_{i,a}^s + e_{i,a}^+ \leq 0, \quad \forall i \in \mathcal{D}, a \in A_i, s \in \mathcal{S} \tag{11.3}$$

$$\zeta_{i,a}^s \geq 0, \quad \forall i \in \mathcal{D}, a \in A_i, s \in \mathcal{S} \tag{11.4}$$

### 3.3 Subtour elimination constraint

The flow conservation and ED/EP energy constraints rule out the infeasible routes that result in the imbalance flows and energy infeasibility. However, those constraints do not prevent the solution which enables an EP to revisit the same ED. This study considers such service loops as subtours. Figure 3(a) shows an example. Specifically, an EP leaves ED(1), conducts a series of local or distant switches to charge several other customers, and eventually, it switches back to serve ED(1) without going back to the depot. Note that this type of subtour is not a route loop on the graph $G(\mathcal{N}, A)$ but a service loop among the EDs. It can be either illegal or legal for the CaaS operation, given that each ED can accept multiple times of services during the trip. Below gives more discussion about these subtours and then develops the subtour elimination constraints for illegal subours in this study.

The Figure 3(b) shows an example of an illegal subtour, where EP(s) serves several EDs from (1) to (4) and return to ED(1) at node $n_1$; its trajectory forms a temporally infeasible service loop. More exactly, EP(s) starts the route by serving ED(1) from nodes $n_1$ to $n_2$. At node $n_2$, it travels to node $n_3$ and serves ED(2) (i.e., an distant switch). After that, EP(s) switches the service to charge ED(3) and then ED(4). Last, the EP(s) arrives at node $n_1$ and switches the service back to ED(1) at $n_1$. This subtour makes EP(s) visit ED(1) at the same node $n_1$ twice but at different time stamps, which is temporally impossible given ED(1) is always on the move toward its destination during this period. Thus, this is an illegal subtour. On the other hand, some subtours are legal in this study. Figure 3(c) shows an example. It starts with charging the ED(1) and moving toward the node $n_1$ and then (local/distant) switches to a series of EDs(2 to 4) and finally switch back to charge ED(1) at $n_2$. This subtour only makes EP(s) revisit ED(1) at node $n_2$ after the first visit at node $n_1$. Given that partial charging is allowed in the CaaS platform, this revisit is permitted and then the subtour is feasible if the corresponding revisits are feasible temporally.

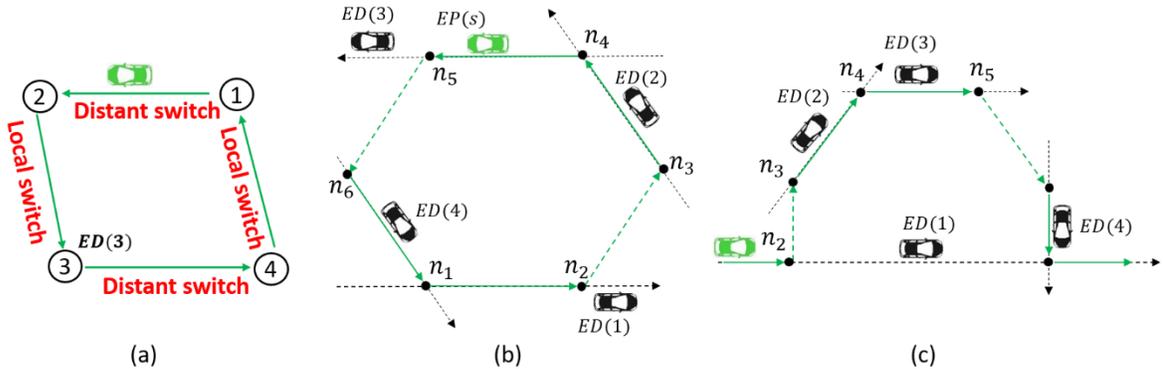

**Figure 3.** (a) Example of the subtour with EDs as nodes and local-/distant-switches as links; (b) An illegal subtour example of (a); (c) A legal subtour example of (a). Green dash-arrow: distant switch.

This study noticed that conventional subtour elimination constraints, e.g., the MTZ constraint (Miller et al., 1960), eliminate all possible revisits since they violate the visiting order constraints. Thus they do not adapt to this study. Motivated by this view, we develop new subtour elimination constraints in (13) and (14). They trace the visiting time of an EP at each node and allow a revisit only if it conducts temporal feasible distant or local switches. Therefore, our subtour elimination constraints only eliminate the illegal subtours but keep the legal ones in this study.

We prove the correctness of our subtour elimination constraints by Theorem 1. To facilitate the development of the proof, we categorize the illegal subtours into three types: the subtour (1) purely formed by distant switches, (2) purely formed by local switches, and (3) mixed with distant and local switches. Then, we prove that constraints (13) (14) can eliminate these three types of illegal subtours (Lemma 1 – 3).

**Lemma 1.** *For any ED(i), a trip satisfying (13) and (14) eliminates the illegal subtours that are only formed by distant switch from/back to ED(i).*

*Proof.* We prove Lemma 1 by showing a contradiction. Let $n$ be the number of EDs involved in the subtour. The proof starts with the case involving two EDs (i.e., $n = 2$), and then extends to a general case with $n \geq 3$. We consider that there exists an illegal subtour purely formed by distant switches as $n = 2$, and this subtour satisfies (13) and (14). For example, $ED(i_1)$ and $ED(i_2)$ in Figure 4(a), in which both $ED(i_1)$ and $ED(i_2)$ travel from node $n_1$ to node $n_2$. Without loss of generality, we assume an EP serves $ED(i_1)$ as they move from node $n_1$ to $n_2$ and conducts distant switch to serve $ED(i_2)$ at node $n_1$ (i.e., $v_{i_1,i_2}^{n_2,n_1,S} = 1$). Then, the EP serves $ED(i_2)$ until it arrives at node $n_2$, from which the EP moves back to node $n_1$ and revisit $ED(i_1)$ (i.e., does a distant switch from $n_2$ to $n_1$ (i.e., $v_{i_2,i_1}^{n_2,n_1,S} = 1$)). Clearly, this subtour is illegal since it makes the EP revisit $ED(i_1)$ at node $n_1$ at different time stamps. We show the conflict resulting from the illegal subtour as follows. Given constraint (14) is satisfied with $v_{i_1,i_2}^{n_2,n_1,S} = v_{i_2,i_1}^{n_2,n_1,S} = 1$, we have (16) and (17).

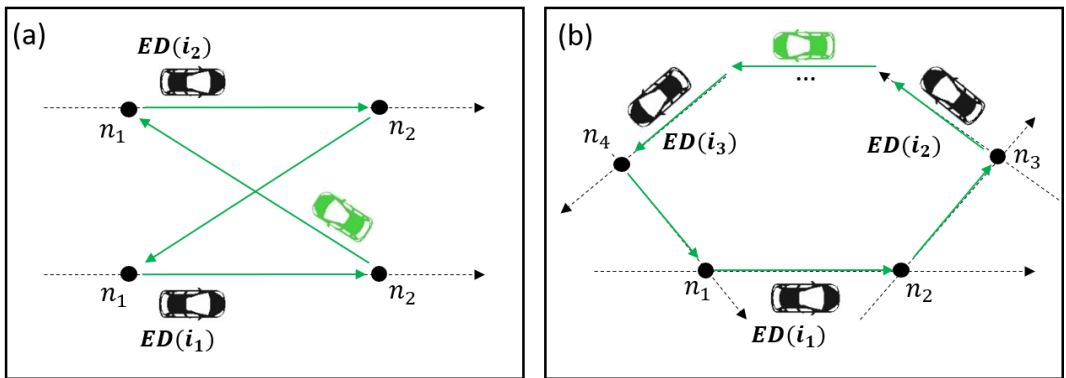

**Figure 4.** Illegal subtours formed by the distant switches only. (a) $n = 2$ and (b) $n \geq 3$.

$$t_{i_2}^{n_1} - t_{i_1}^{n_2} - t_{n_1 \to n_2} \geq 0, \tag{16}$$

$$t_{i_1}^{n_1} - t_{i_2}^{n_2} - t_{n_2 \to n_1} \geq 0. \tag{17}$$

Adding (16) and (17) together and restructuring the result, we have

$$-2t_{n_2 \to n_1} \geq \left(t_{i_1}^{n_2} - t_{i_1}^{n_1}\right) + \left(t_{i_2}^{n_2} - t_{i_2}^{n_1}\right) > 0. \tag{18}$$

Equation (18) indicates that $t_{n_2 \to n_1} < 0$, which contradicts with the fact $t_{n_2 \to n_1} > 0$. Therefore, such subtour does not exist.

When the case with $n \geq 3$ (see Figure 4(b) for an example), we assume an EP serves ED($i_1$) from node $n_1$ to $n_2$ and conducts distant switch to serve ED($i_2$) ($v_{i_1,i_2}^{n_2,n_3,S} = 1$). Then, the EP continues to serve several other EDs by multiple distant switches and finally revisit ED($i_1$) at node $n_1$ after it leaves ED($i_3$) ($v_{i_3,i_1}^{n_4,n_1,S} = 1$) from node $n_4$. Focusing on the feasibility of the first and last distant switches associated with ED($i_1$), we have (19)(20) below given constraint (14) is satisfied with $v_{i_1,i_2}^{n_2,n_3,S} = v_{i_3,i_1}^{n_4,n_1,S} = 1$.

$$t_{i_2}^{n_3} - t_{i_1}^{n_2} - t_{n_2 \to n_3} \geq 0 \tag{19}$$

$$t_{i_1}^{n_1} - t_{i_3}^{n_4} - t_{n_4 \to n_1} \geq 0 \tag{20}$$

Adding (19) and (20) together and restructuring the result, we have (21) below.

$$-t_{n_4 \to n_1} - t_{n_2 \to n_3} \geq \left(t_{i_1}^{n_2} - t_{i_1}^{n_1}\right) + \left(t_{i_3}^{n_4} - t_{i_2}^{n_3}\right) > 0 \tag{21}$$

(21) indicates that $t_{n_4 \to n_1} + t_{n_2 \to n_3} < 0$, which contradicts with $t_{n_2 \to n_3} > 0, t_{n_4 \to n_1} > 0$. Therefore, the subtour does not exist. We thus complete the proof.

**Lemma 2.** *For any ED(i), a trip satisfying (13) and (14) eliminates subtours that are only formed by local switch from/back to ED(i).*

*Proof.* The procedure to prove this lemma is similar to the proof of Lemma 1. We leave it in Appendix D.

**Lemma 3.** *For any ED(i), a trip satisfying (13) and (14) eliminates subtours that are formed by distant switch from/back to and local switch back to/from ED(i).*

*Proof.* The procedure to prove this lemma is similar to the proof of Lemma 1. We leave it in Appendix D.

**Theorem 1.** *Constraints (13) (14) are valid subtour elimination constraints for the mE2-VRP model.*

*Proof.* Theorem 1 is immediately proved by Lemma 1-3.

## 4. Clustering-Aided Clarke and Wright Savings

The mE2-VRP is a mixed-integer program subject to many constraints, which is an NP-hard problem in general. Moreover, the mE2-VRP considers on-the-move electricity delivery and allows an EP to revisit the same ED multiple times. These unique features introduce extra complexity to the mE2-VRP as compared to the similar problems in the literature such as electric vehicle dial-a-ride problems (e-DARP) (Bongiovanni et al., 2019), and electric (Schneider et al., 2014) and green VRP (Erdoğan and Miller-Hooks, 2012). Finding the optimal solution of the mE2-VRP may be a computationally impractical task. This study, therefore, focuses on efficiently searching a good-quality solution to satisfy the timely service requirement of the CaaS platform. Accordingly, we design the Clustering-aided Clarke and Wright Savings (CCWS) algorithm as a heuristic approach customized towards the complications of the mE2-VRP. The main idea of the CCWS is introduced as follows.

First of all, it is noticed that the computation complexity of the mE2-VRP increases significantly as more EDs are involved in the problem. This is because those EDs will introduce more nodes and arcs to the graph $G(\mathcal{N}, A)$ and complicates its topology. As a result, it makes the feasible trips of EPs (feasible solution) enormous. However, we also can observe that some EDs, constrained by their trip plans, will not have the potential to be served by the same EP. Then, they can be treated separately. By understanding this feature in depths, the CCWS first strategically splits all EDs into multiple clusters. It enables us to decompose the master problem (i.e., mE2-VRP) into a number of subproblems (i.e., sub-mE2-VRPs). Each subproblem has a small instance size and can be efficiently solved by existing commercial solvers such as Gurobi. Next, we solve each sub-mE2-VRP simultaneously by parallel computing. The tours of the EPs in each cluster are considered as the seed tours. All seed tours together form the first feasible solution of the master mE2-VRP. Figure 5 (a) shows an example. Last, this study further improves this initial feasible solution by strategically merging the seed tours over all clusters. It will help reduce the size of the EP fleet and represents a better solution for the master mE2-VRP (see the illustration of the idea in Figure 5 (b)). In summary, the CCWS algorithm includes three key steps as follows.

**Step 1:** Form the ED clusters and decompose the master problem into the sub-mE2-VRPs.

**Step 2:** Parallelly solve each sub-mE2-VRP and generate seed tours.

**Step 3:** Merge seed tours to explore a better solution.

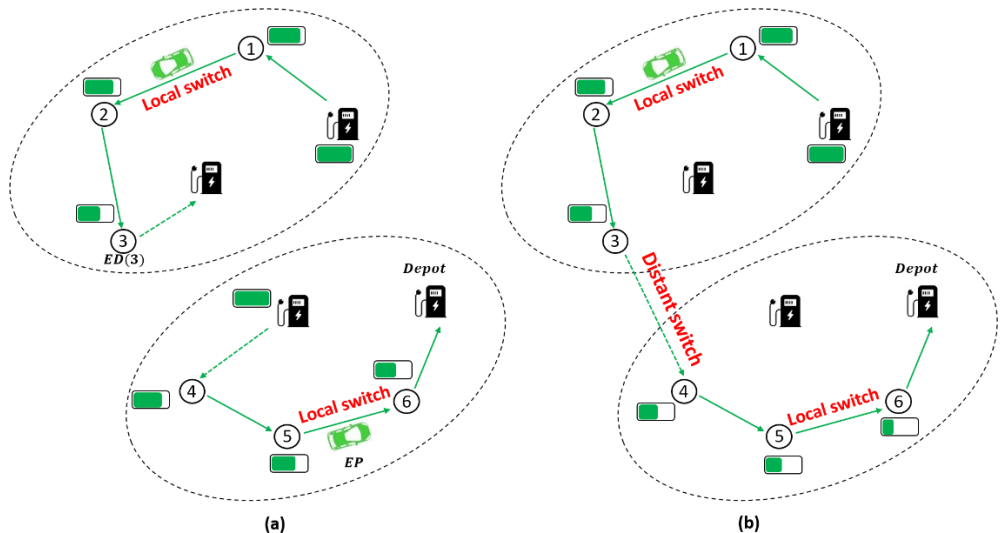

**Figure 5.** Schematic representation of the CCWS algorithm.

There are several critical technical issues in each of the three steps above, which will affect the performance of the CCWS significantly. Mainly, to ensure a good initial solution formed by the seed tours, the ED clusters should be strategically formed so that we can minimize the compromise of the global optimality resulting from the decomposition. Moreover, the CCWS requires careful efforts to solve the subproblems and develop an appropriate merging algorithm to secure the feasibility and efficiency of the improved solution. We discuss these technical challenges and our approaches in detail in the following sections.

**4.1 Forming ED clusters**

We first discuss how to form the ED clusters so that they can facilitate the solution exploration as much as possible. It is noticed that distant switches are expensive since each distant switch causes extra idle energy and time loss. Therefore, an optimal solution prefers local switches rather than distant switches. Following this thought, we recognized that if we cluster the EDs with high potential of local switches and then explore the seed tours in each cluster to form the initial solution, this decomposition is likely to lead to a good feasible solution without over sacrificing global optimality. Invoked by these observations, this study develops a clustering algorithm to form the ED groups, factoring in the possibility of the local switch occurring between two EDs according to their trip plans. More exactly, the clustering algorithm seeks to put two EDs in one cluster if their trips present a high chance to enable potential local switches; otherwise, they will be separated into different clusters. Note that we temporally ignore the limit of the energy inventory as we form the ED clusters.

Our clustering algorithm consists of two parts. (1) Measure the local switch opportunities between any pair of the EDs considering their routes and trip schedules. (2) Forming the ED clusters with high local switch potential among the EDs in each cluster. However, to the best of our knowledge, there is not an existing approach to quantify such local switch potential between a pair of EDs. It is uniquely introduced by this study. Therefore, we need to develop our own measurement, and then strategically cluster the EDs upon it.

To quantify the chance of a local switch, we consider two EDs such as $i$ and $j$, both of which will pass node $n$ in $G(\mathcal{N}, A)$. If ED($i$) under the service of an EP arrives at node $n$ is earlier than ED($j$) does (i.e., $t_n^j - t_n^i \geq 0$), then potentially the EP can locally switch from ED($i$) to serve ED($j$) at node n. All such ED($j$) form the Local-Switch Candidate (LSC) set for ED($i$) at node $n$, which is labeled as $J_{i,n}$. The size of this set can help us to measure the opportunities that local switches will occur between a given ED and all other EDs. An ED with a larger such set will have more opportunities to enable an EP on duty with it to do a local switch and continue to serve the other EDs at a node. Accordingly, we formally define the LSC set of ED(i) at node $n$ and the LSC set for all EDs in consideration in (22).

$$J_{i,n} = \{ED(j)|t_n^j - t_n^i \geq 0, j \in D_n, j \neq i\}, \quad i \in \mathcal{D}, n \in \mathcal{N}_i \setminus i_0; J = \{J_{i,n}\}_{\{i \in \mathcal{D}, n \in \mathcal{N}_i \setminus i_0\}} \tag{22}$$

Recall that $t_n^i = \tilde{t}_0^i + t_{i_0 \to n} + \tau_i$, $t_n^j = \tilde{t}_0^j + t_{j_0 \to n} + \tau_j$ in (15), and $\tau_i$ and $\tau_j$ are continuous variable respectively representing the waiting time of ED(i) and ED(j), $0 \leq \tau_i \leq \bar{\tau}_i, 0 \leq \tau_j \leq \bar{\tau}_j$. Given $\tau_i$ and $\tau_j$ are unknown decision variables in the mE2-VRPs mathematical model, we are not able to directly apply set $J_{i,n}$ in the development of the clustering algorithm. Thus, we consider another equivalent set $\bar{J}_{i,n}$ defined by (23) which explicitly formulates the LSC set of an ED at a node, using known parameters $\tilde{t}_0^i$, $\tilde{t}_0^j$, $t_{i_0 \to n}$ and $t_{j_0 \to n}$. Lemma 4 proves the equivalence of $\bar{J}_{i,n}$ and $J_{i,n}$ as well as $J$ and $\bar{J}$.

$$\bar{J}_{i,n} = \{ED(j)|(\tilde{t}_0^j + t_{j_0 \to n}) - (\tilde{t}_0^i + t_{i_0 \to n}) \geq -\bar{\tau}_j, j \in D_n, j \neq i\}, i \in \mathcal{D}, n \in \mathcal{N}_i \setminus i_0; \tag{23}$$
$$\bar{J} = \{\bar{J}_{i,n}\}_{\{i \in \mathcal{D}, n \in \mathcal{N}_i \setminus i_0\}}$$

**Lemma 4.** *The two sets $\bar{J} = \{\bar{J}_{i,n}\}_{\{i \in \mathcal{D}, n \in \mathcal{N}_i \setminus i_0\}}$ and $J = \{J_{i,n}\}_{\{i \in \mathcal{D}, n \in \mathcal{N}_i \setminus i_0\}}$ are equivalent.*

*Proof.* We first prove $J \subseteq \bar{J}$. To do it, our proof shows that any ED($j$) in the set of $J_{i,n}$ belongs to $\bar{J}_{i,n}$ too, i.e., $J_{i,n} \subseteq \bar{J}_{i,n}, \forall i \in \mathcal{D}, n \in \mathcal{N}_i \setminus i_0$. According to (22) and (15), we have (24)-(26) for any $ED(j) \in J_{i,n}$.

$$t_n^i - t_n^j \geq 0 \Leftrightarrow (\tilde{t}_0^j + t_{j_0 \to n} + \tau_j) - (\tilde{t}_0^i + t_{i_0 \to n} + \tau_i) \geq 0 \tag{24}$$
$$0 \leq \tau_i \leq \bar{\tau}_i \tag{25}$$
$$0 \leq \tau_j \leq \bar{\tau}_j \tag{26}$$

According to (25) and (26), we restructure (24) to get (27) below for $ED(j) \in J_{i,n}$.

$$(\tilde{t}_0^j + t_{j_0 \to n}) - (\tilde{t}_0^i + t_{i_0 \to n}) \geq \tau_i - \tau_j \geq -\bar{\tau}_j \tag{27}$$

Then, (23) and (27) together indicate that $ED(j) \in \bar{J}_{i,n}$ for any $ED(j) \in J_{i,n}$. Therefore, we have the conclusion in (28).

$$J_{i,n} \subseteq \bar{J}_{i,n}, \forall i \in \mathcal{D}, n \in \mathcal{N}_i \setminus i_0 \Leftrightarrow J \subseteq \bar{J} \tag{28}$$

We next prove that $\bar{J} \subseteq J$ by showing that any $ED(j)$ in the set of $\bar{J}_{i,n}$ belongs to $J_{i,n}$ too, i.e., $\bar{J}_{i,n} \subseteq J_{i,n}, \forall i \in \mathcal{D}, n \in \mathcal{N}_i \setminus i_0$. To do that, we denote a set for all possible waiting time tuples $(\tau_i, \tau_j)$ for $ED(j)$ in $J_{i,n}$ as

$$\Omega_{i,j}^n = \{(\tau_i, \tau_j) | (\tilde{t}_0^j + t_{j_0 \to n} + \tau_j) - (\tilde{t}_0^i + t_{i_0 \to n} + \tau_i) \geq 0, 0 \leq \tau_i \leq \bar{\tau}_i, 0 \leq \tau_j \leq \bar{\tau}_j\}, i \in \mathcal{D}, n \in \mathcal{N}_i \setminus i_0, j \in D_n, j \neq i$$

Clearly, if we can find a pair of waiting time $(\tau_i, \tau_j)$ for $ED(i)$ and $ED(j)$ in $\Omega_{i,j}^n$, then $ED(j)$ is in the set of $J_{i,n}$. Considering any $j \in \bar{J}_{i,n}$, we have (29) according to (23).

$$(\tilde{t}_0^j + t_{j_0 \to n} + \bar{\tau}_j) - (\tilde{t}_0^i + t_{i_0 \to n}) \geq 0 \tag{29}$$

(29) indicates that there exist a waiting time tuple: $\tau_j = \bar{\tau}_j$ and $\tau_i = 0$ such that $(\tau_j, \tau_i) \in \Omega_{i,j}^n$. Therefore, we conclude that any $j \in \bar{J}_{i,n}$ satisfies $j \in J_{i,n}$, and thus we have the result in (30).

$$\bar{J}_{i,n} \subseteq J_{i,n}, \forall i \in \mathcal{D}, n \in \mathcal{N}_i \setminus i_0 \Leftrightarrow \bar{J} \subseteq J \tag{30}$$

Combining the results in (28) and (30), we conclude Lemma 4 #.

The set $\bar{J}_{i,n}$ help us understand the chance that a local switch may occur between $ED(i)$ and all other EDs at a single node $n$. It is not sufficient for the design of the clustering algorithm since it still cannot capture the opportunities for an EP to perform local switches between a pair of EDs. To solve this issue, we develop Equation (31) based upon the set $\bar{J}$, which measures the local-switch opportunity between a pair of EDs such as $i$ and $j$ by the number of potential local switches from ED(i) to ED(j) along the route of ED(i).

$$p(i,j) = \begin{cases} \sum_{n \in \mathcal{N}_i \setminus i_0} \mathbb{1}\{j \in J_{i,n}\}, \forall i,j \in D, i \neq j \\ 0, \quad i = j \end{cases} \text{ and } P = [p(i,j)]_{i,j \in D} \tag{31}$$

where $\mathbb{1}\{j \in J_{i,n}\}$ is an indicator function, takes value 1 when $j \in J_{i,n}$.

Last, built upon the measurement of the local-switch opportunity, we design the clustering algorithm to form ED clusters. Mainly, the clustering algorithm starts with the ED pair that has the highest local switch opportunities and puts them in one cluster, then iteratively clustering the next promising ED in the same cluster until no more ED can be locally switched to. Appendix B presents and explains the steps of the pseudocode of our clustering algorithm.

## 4.2 Generating seed tours to form a feasible solution

We next discuss the second step of the CCWS, which seeks to solve the sub-mE2-VRP for each cluster efficiently. Our experiments noticed that each subproblem may still involve many EDs (such as >100 EDs) as we implement the CaaS in a large network. This is because the clustering algorithm is an unsupervised learning algorithm and cannot pre-determine the size of each ED cluster. As a result, the corresponding sub-mE2-VRP may still has many variables and constraints and is difficult to be solved within planning interval length $T$. To address this challenge, the study develops the scheme shown in (32). It disables distant switches for each cluster $c$ to sustain the computation efficiency and solution quality by leveraging two critical features of the problem. First of all, we observed that distant switch services contribute most to the model complexity since each of them introduces more variables and constraints than a local switch does in the mathematical model. Next, considering the clusters are formed with particular interest to promote local switch services in each cluster, we have a strong reason to believe that an optimal solution of each subproblem will involve more local switches than distant switch services. Thanks to these merits of the ED clusters, disabling distant switches in a subproblem will benefit the computation load without over sacrificing the solution quality[3] for solving the sub-mE2-VRP. In the meantime, the complexity of the sub-mE2-VRP($c$) is significantly reduced and these subproblems can be quickly solved (find seed tours) by commercial solvers using parallel computation.

**sub-mE2-VRP($c$):** $\{\min F(O), s.t. \ v_{i,j}^{n,m,s} = 0, \forall v \in V; s \in \mathcal{S}, i \in \mathcal{D}_c; \ (2) \sim (13)\}$ (32)

Putting all the seed tours together we can quickly obtain the initial feasible solution of the mE2-VRP. Mathematically, we denote $\mathcal{S}_0 \subseteq \mathcal{S}$ as the set of the EPs that are dispatched in the initial feasible solution. Then, each EP($s$) in $\mathcal{S}_0$ is associated with a seed tour $H(s), s \in \mathcal{S}_0$ (see Figure 6). Accordingly, Equation (33) includes all necessary information of a seed tour $H(s)$ required for merging which will be discussed later, where $n$ and $t_0^{*s}$ represents the locations and the time that EP($s$) start the service in the seed tour. Here, we use (*) to represent the optimal solution of sub-mE2-VRP. Similarly, $m$ and $t_e^{*s}$ represents the locations and the time that EP(s) completes the service with the last ED in the seed tour.

$$H(s) = (n, t_0^{*s}, m, t_e^{*s}, E_s^*)$$ (33)

Last, $E_s^*$ represents the total energy loss of EP(s) during this seed tour. Note that $E_s^*$ excludes the idle energy loss from a depot to the first ED, and from the last ED to the depot.

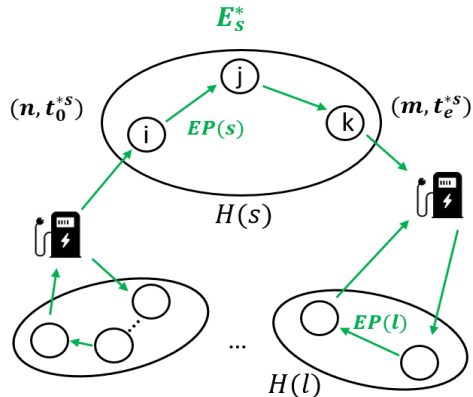

**Figure 6.** An example of seed tours.

---

[3] Even though we disable the distant switch services in the master mE2-VRP, it is still a large-scale mixed integer programming and is very hard to be solved. In addition, without the merits of the clustering algorithm, this simplification will comprise the optimality significantly.

Now consider any EP(s) in $\mathcal{S}_0$, we denote $\mathcal{D}_s^* \subseteq \mathcal{D}$ as the set of the EDs served by the EP(s), e.g., $\mathcal{D}_s^* = \{i, j, k\}$ in an example of Figure 6. Specifically, we use $s(n)$ to represent the $n$-th ED that EP(s) serves in the seed tour, e.g., $s(1) = i, s(2) = j$, and the last ED that EP(s) serves is denoted as $s(-1)$, e.g. $s(-1) = k$. With these notations, we introduce how to merge the seed tours and improve the initial feasible solution in the next section.

### 4.3 Merging seed tours to refine the feasible solution

The initial feasible solution obtained in Step 2 excludes any distant switch service between the EP trips formed in individual clusters. Given the mE2-VRP seeks to fulfill the service requests with a minimum fleet size, this study next wants to further improve this feasible solution by merging those seed tours by adding back the distant switches. To present our main idea, let's consider two seed tours, $H(s) = (n_1, t_0^{*s}, n_2, t_e^{*s}, E_s^*)$ and $H(l) = (m_1, t_0^{*l}, m_2, t_e^{*l}, E_l^*)$, $s, l \in \mathcal{S}_0$. The merging algorithm seeks to have an EP serve the EDs in the seed tour $H(s)$ first and then $H(l)$ without violating the energy feasibility (Rule 1) and delay tolerance (Rule 2). Each successful merging will reduce the fleet size of the initial feasible solution by one. More exactly, we would like to merge two seed tours to one if an EP has enough energy to fulfill the service requests of the EDs in the other seed tour such as $H(l)$ after it finishes with the seed tour on duty such as $H(s)$. Clearly, each merging introduces one more distant switch (and also delay), by which the EP travels from the tail of $H(s)$ to the head of $H(l)$. Therefore, a successful switch needs to also satisfy the delay tolerance required by the EDs in the second seed tour of $H(l)$. Following this idea, we developed the merging algorithm according to the two merging rules in (34) and (35).

**Rule 1:** $e_0^S - (e_{s_0,n_1}^- + E_s^*) - e_{n_2,m_1}^- - (E_l^* + e_{m_2,p_{m_2}}^-) \geq e_s$ (34)

**Rule 2:** $(t_e^{*s} + t_{n_2 \to m_1}) - t_0^{*l} \leq \min\{\bar{\tau}_i - \tau_i^*, \forall i \in \mathcal{D}_l^*\}$ (35)

**Rule 1** ensures that EP(s) has enough energy to fulfill energy requests of the EDs in $H(l)$. Mainly, the first item is the initial energy inventory of EP(s), the second and fourth item are the energy losses associated with the seed tours $H(s)$ and $H(l)$, respectively. The third item is the energy required to do the distant switch (merging). After serving all the EDs in both seed tours, the EP should have enough energy left to go back to depot. **Rule 2** guarantees that EP(s) can timely meet every EDs in the second seed tour without violating their waiting tolerance. The left-hand-side of Rule 2 is the extra delay resulting from the merging, which involves the difference of two items. The first one is the time that EP(s) can start serving the head of $H(l)$ if these two seed tours, $H(s)$ and $H(l)$ are merged, and the second item is the time that the head of $H(l)$ is served by EP(l) in the pre-merging solution. The right-hand-side of Rule 2 is the upper bound of the extra delay imposed by the EDs in $H(l)$. Each ED(i) in $H(l)$ has a tolerable delay $\bar{\tau}_i$ and an actual delay, $\tau_i^*$ in the pre-merging solution. Therefore, $\bar{\tau}_i - \tau_i^*$ is the extra delay that ED(i) can tolerate if two seed tours are merged. The minimum of these extra delays bounds the extra delay from merging so that the waiting tolerance of all the EDs in $H(l)$ will not be violated.

The algorithm is interested in merging as many as possible seed tours to one EP trip so that we can minimize the size of the EP fleet. To achieve this, we notice that the departure time after an EP finishes with the last ED in a merged tour is critical. To facilitate our discussion, we denote a merged tour in the example above as $\bar{H}(s-l)$ and this departure time from the tour as $t_e^{s-l}$. We found that an EP leaving the $\bar{H}(s-l)$ earlier will have more opportunities to accommodate another seed tour following the tour of $\bar{H}(s-l)$ in one trip. Consequently, our merging algorithm gives high priority to merge a pair of seed tours which leads to a smaller value of $t_e^{s-l}$. Appendix C provides the pseudocode of the merging algorithm. Here, we briefly discuss the main steps of our merging algorithm. We first find all feasible seed tours for merging according to Rule 1 and Rule 2 and rank their priority. Following that, we select the best feasible seed tours to merge according to their priority. Last, considering the merged tour as a new tour, repeat the first two steps until no more feasible merging pair can be found.

### 4.4 The feasibility of the CCWS algorithm

This section discusses if the CCWS algorithm can guarantee a feasible solution to the mE2-VRP model. Clearly, the initial solution formed by gathering the seed tours from solving the subproblems is feasible to the mE2-VRP by putting all distant switch variables equal to zeros (i.e., making constraint (14) always correct). Our doubt is raised by the merging process which adds back distant switches according to constraints (34) and (35). More exactly, constraint (34) used to ensure the energy feasibility in the CCWS is different from the energy inventory constraint (12) in the mE2-VRP. Similarly, constraint (35) used to secure the delay tolerance in the CCWS is also different from constraints (13) and (14), which ensure the feasible local and distant switches according to the EDs' delay tolerance. In view of these issues, this section proves that a CCWS solution satisfying constraints (34) and (35) will also satisfy constraints (12), (13), and (14), respectively. In this way, we can guarantee that the CCWS algorithm finds a good feasible solution to the mE2-VRP. Lemma 5 below summarizes this thought.

**Lemma 5.** *An EP tour obtained by merging the seed tours of the CCWS algorithm according to constraints (34) and (35) will satisfy constraints (12), (13), and (14) in the mE2-VRP.*

*Proof.* First, the seed tours of the CCWS are obtained by solving the sub-mE2-VRPs. Therefore, they satisfy constraints (12), (13), and (14) in the mE2-VRP. Next, we further prove that merging these seed tours according to (34) and (35) won't violate constraints (12), (13), and (14).

We first prove that the EP tour obtained satisfying constraint (34) won't violate constraint (12). Clearly, by taking $E_s = e^-_{s_0 \to n_1} + E^*_s + e^-_{n_2 \to m_1} + E^*_l + e^-_{m_2 \to p_{m_2}}$ (the left side of (34)), constraint (34) indicates that $e^S_0 - E_s \geq \underline{e}_s$, which satisfies constraint (12).

We next prove a tour that satisfies constraint (35) also satisfies constraints (14). To conduct this proof, we consider two seed tours, such as $H(s) = (n_1, t^{*s}_0, n_2, t^{*s}_e, E^*_s)$ and $H(l) = (m_1, t^{*l}_0, m_2, t^{*l}_e, E^*_l)$. The CCWS algorithm merges them according to the constraint (35). We denote the last ED in tour $H(s)$ by ED(i) and the first ED in $H(l)$ by ED(j)) (i.e., $s(-1) = i$, $l(0) = j$). They are respectively associated with the waiting time $\tau^*_i \in [0, \bar{\tau}_i]$ and $\tau^*_j \in [0, \bar{\tau}_j]$. The proof first shows that the distant switch from ED(i) to ED(j) (i.e., $v^{n_2,m_1,s}_{i,j} = 1$) satisfies constraint (14). Combining with (15), this is equivalent to prove (36) holds. Namely, there exists feasible waiting delays satisfying the delay requests of ED(i) in $H(s)$ and ED(j) in $H(l)$ such that ED(j) can receive timely service when EP(s) conduct a distant switch from ED(i) to ED(j).

$$\exists \tau_i \in [0, \bar{\tau}_i], \tau_j \in [0, \bar{\tau}_j], \ s.t \ \left(\tilde{t}^j_0 + t_{j_0 \to m_1} + \tau_j\right) - \left(\tilde{t}^i_0 + t_{i_0 \to n_2} + \tau_i\right) - t_{n_2, m_1} \geq 0 \quad (36)$$

To prove the correctness of (36), we consider (35) is satisfied as we merge $H(s)$ and $H(l)$. Given $t^{*s}_e = t^i_{n_2}$ and $t^{*l}_0 = t^j_{m_1}$, we plug (15) into (35), and obtain (37).

$$\left(\tilde{t}^j_0 + t_{j_0 \to m_1} + \tau^*_j\right) - \left(\tilde{t}^i_0 + t_{i_0 \to n_2} + \tau^*_i\right) - t_{n_2 \to m_1} + \min\{\bar{\tau}_i - \tau^*_i, \forall i \in \mathcal{D}^*_l\} \geq 0 \quad (37)$$

Now, we can take $\tau_i = \tau^*_i \in [0, \bar{\tau}_i]$, and $\tau_j = \tau^*_j + \min\{\bar{\tau}_i - \tau^*_i, \forall i \in \mathcal{D}^*_l\} \leq \tau^*_j + \bar{\tau}_j - \tau^*_j = \bar{\tau}_j$ and make (36) hold. Therefore, the merging satisfying (35) won't violate the constraint (14).

Last, we prove a tour that satisfies constraint (35) will make constraints (13) hold too. Clearly, the merging algorithm won't affect the trip plan of the first seed tour such as $H(s)$ in our example. And then all the local switches in $H(s)$ satisfy constraint (13). However, we need to show that the merging can also make the local switches in the second seed tour such as $H(l)$ satisfy constraint (13). Here, we only prove that the first local switch in the seed tour $H(l)$ is feasible. Follow the same logic, the feasibility of other local switches in the $H(l)$ can be easily proved.

We denote the first and the second EDs in the seed tour $H(l)$ by ED($j$) and ED($k$). The EP($s$) conducts a local switch (the first local switch in $H(l)$) at node $c$ from the ED($j$) to the ED($k$), i.e., $u_{j,k}^{c,s} = 1$. We show this local switch is feasible to constraint (13). Combining with (15), this is equivalent to prove (38) holds. Namely, there exist feasible waiting delays satisfying the delay requests of ED(j) and ED(k) in $H(l)$, such that when two seed tours merge, ED(k) can be timely served by EP(s).

$$\exists \tau_j \in [0, \bar{\tau}_j], \tau_k \in [0, \bar{\tau}_k], \ s.t \ \left(\tilde{t}_0^k + t_{k_0 \to m_3} + \tau_k\right) - \left(\tilde{t}_0^j + t_{j_0 \to m_3} + \tau_j\right) \geq 0 \tag{38}$$

We have known that the initial solutions $\tau_k^*$ and $\tau_j^*$ satisfy constraint (13). By (15), constraint (13) is equivalent to (39).

$$t_{m_3}^k - t_{m_3}^j = \left(\tilde{t}_0^k + t_{k_0 \to m_3} + \tau_k^*\right) - \left(\tilde{t}_0^j + t_{j_0 \to m_3} + \tau_j^*\right) \geq 0 \tag{39}$$

Now let's take $\tau_j = \tau_j^* + \min\{\bar{\tau}_i - \tau_i^*, \forall i \in \mathcal{D}_l^*\} \leq \bar{\tau}_j$ and $\tau_k = \tau_k^* + \min\{\bar{\tau}_i - \tau_i^*, \forall i \in \mathcal{D}_l^*\} \leq \tau_k^* + \bar{\tau}_k - \tau_k^* = \bar{\tau}_k$. The constraint in (38) becomes $\left(\tilde{t}_0^k + t_{k_0 \to m_3} + \tau_k^*\right) - \left(\tilde{t}_0^j + t_{j_0 \to m_3} + \tau_j^*\right) \geq 0$, which is satisfied according to (39). Therefore, there exists $\tau_j \in [0, \bar{\tau}_j], \tau_k \in [0, \bar{\tau}_k]$, such that (38) holds. Now, we can conclude that the local switches after merging are still feasible to constraint (13). Wrapping up all conclusions above, we complete this proof. It also indicates that the merging won't introduce infeasibility.

**Theorem 2.** *The CCWS algorithm finds a feasible solution to the mE2-VRP.*

*Proof.* First, the initial solution obtained from the CCWS algorithm is a collection of the solutions of the sub-mE2-VRPs defined by (32), whose constraints are subject to the mE2-VRP. Therefore, this initial solution is feasible to the mE2-VRP. Then, Lemma 5 proves that the merging instrument for involving the distant switches to improve the quality of the solution without violating energy inventory constraint in (12), and delay tolerance (subtour-elimination) constrains in (13) and (14) in the mE2-VRP. Therefore, the CCWS solution is a feasible solution to the master mE2-VRP.

## 5. Numerical study

This section sets up numerical studies to validate the performance of the model and algorithm developed in this study. Built upon that, we further explore insights and the applicability of the CaaS platform in citywide and statewide implementations. To do that, we conducted two sets of numerical experiments, which respectively built upon the Chicago sketch network (Mayer, 1961) and the Florida statewide network. All numerical experiments are run on a DELL Precision 3630 Tower with 3.60GHz of Intel Core i9-9900k CPU, 8 cores and 16 GB RAM in a Windows environment.

The two sets of experiments share the default setting for some parameters of the EDs and the EPs. Specifically, the ED battery capacity is assumed to be 90 kWh, with consumption rate $\varpi = 0.4$ kWh/mile (max driving range of 225 miles), which corresponds to the battery performance of Tesla model S. We consider that the mE2 charging service will relax the range anxiety, so that people would feel comfortable to leave home or workplaces without having their cars fully charged. Therefore, the initial battery inventory of the EDs is randomly generated above the safety inventory (2 kWh), and below the energy required to complete the trip, so that each ED in this experiment needs to have at least one recharging service. The EP fleet has the same battery parameters as the EDs, except that each EP can carry more than one battery. The experiments have each EP equipped with two full batteries as the default setting, unless otherwise stated. We further assume a power transfer rate $\eta$ of 55 kW, and transfer efficiency of 0.9 in the mE2 platform (Abdolmaleki et al., 2019; Chakraborty et al., 2020).

### 5.1 Numerical study 1: Chicago sketch network
### 5.1.1 Experiment settings
We first conduct the experiments built on the Chicago sketch network. The objective of the experiments is three-fold: (i) evaluate the computation performance of the CCWS algorithm, (ii) examine the traffic impact of the CaaS platform at the network level, and (iii) analyze the CaaS platform performance under different input parameters. The Chicago sketch network as shown in Figure 7 consists of 933 nodes and 2,950 links. The network is developed by the Chicago Area Transportation Study (CATS) (Mayer, 1961). It is a realistic yet aggregated representation of the Chicago region. We use it as a testbed to evaluate the CaaS platform over citywide demands.

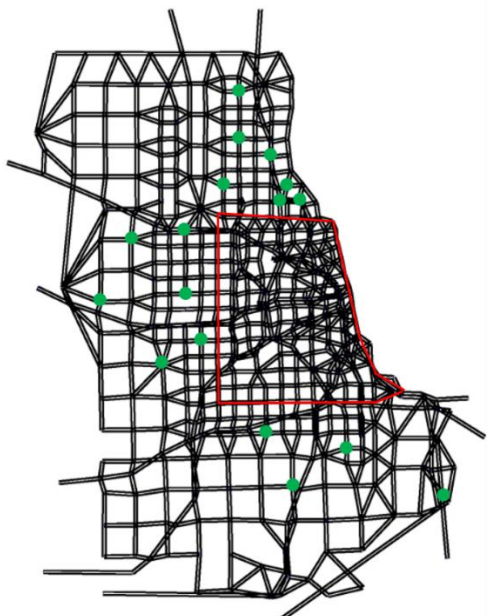

**Figure 7.** Sioux Fall network and the EP depots (green nodes).

More exactly, the experiments simulate the ED trips from a field Chicago O-D trip dataset (Stabler et al., 2018), which consists of 1,260,910 trips. The average ED trip length is around 11 miles. We select the EP recharging depots (green nodes in Figure 7) according to the Nissan dealers locations in Chicago, because the dealers that sell electric vehicles usually have available chargers (Chen and Nie, 2015). Note that this study avoids the CBD area since it usually experiences heavy traffic and not proper to setup the CaaS service.

### 5.1.2 Computation performance of the CCWS algorithm
We first examine the computation performance of the CCWS algorithm under different instances with the size (the number of ED trips) varying from 60 to 10,000. The experiment results of all instances are summarized in Table 1, in which the performance of the CCWS heuristic algorithm is compared with an existing commercial solver, Gurobi, in terms of the solution time and the EP fleet size (the objective function value; the smaller the better). Note that Gurobi is not efficient for the instance with large size. Therefore, we set up solution time limits as an hour and two hours in Gurobi.

**Table 1:** Computation performance of CCWS algorithm under the default setting

| Instance size | Gurobi | | CCWS | | Service rate |
|---|---|---|---|---|---|
| | Solution (1h) | Solution (2h) | Solution | Solution time/s | |
| 60 | 14 | 14 | 14 | 6.02 | 4.29 |
| 80 | 20 | 19 | 20 | 6.23 | 4.00 |
| **100** | **26** | **24** | **24** | **7.20** | **4.17** |
| 120 | - | - | 24 | 8.50 | 5.00 |
| 140 | - | - | 26 | 9.52 | 5.38 |
| 200 | - | - | 32 | 12.59 | 6.25 |
| 300 | - | - | 47 | 17.50 | 6.38 |
| 600 | - | - | 88 | 30.66 | 6.82 |
| 1,200 | - | - | 182 | 54.91 | 6.59 |
| 2,000 | - | - | 289 | 107.50 | 6.92 |
| 4,000 | - | - | 589 | 196.78 | 6.79 |
| 6,000 | - | - | 948 | 374.86 | 6.33 |
| 8,000 | - | - | 1,211 | 536.93 | 6.61 |
| 10,000 | - | - | 1,524 | 691.64 | 6.56 |

The results in Table 1 indicate that the Gurobi solver can only find a solution for an instance with the size less than 100 EDs within 2 hours, but fails as the size of an instance is larger than 100 EDs. Thus, we cannot rely on this commercial software to implement the CaaS platform in an urban area. On the contrary, the CCWS algorithm finds comparable solutions for all the instances up to 10,000 with a more competitive computation efficiency. For example, the CCWS takes only a few seconds to finds a solution for the instance with 100 EDs, and less than 15 minutes for the instance with 10,000 EDs. Moreover, the CCWS is able to find a solution with comparable optimality to that found by Gurobi within two hours for instances less than 100, and even better as compared with Gurobi's one-hour solution (e.g., 24 EPs by the CCWS vs. 26 EPs by Gurobi solver for 100 EDs instance in Table 1). Therefore, the experiments confirmed that the CCWS algorithm addresses the scalability issue of the mE2-VRP model. This super capability makes the CaaS platform possible for an online application. More importantly, it enables this study to set up citywide or statewide experiments to evaluate the system performance of the CaaS platform with realistic settings under different traffic circumstances.

Beyond the fleet size, the experiments also record and analyze the EP service rate under different instances (the last column in Table 1), which is the ratio of the number of EDs to the EP fleet size. This service rate represents the number of EDs that one EP can serve on average. Thus it reflects the service efficiency and implies the potential traffic impact. We will discuss this result together with other experiment results in the next section 5.1.3.

### 5.1.3 Traffic impact analysis

One concern for implementing the CaaS platform is the extra traffic introduced by the EP fleet. Table 1 demonstrates that the EP service rate increases with the instance size, and gets a plateau when the instance size is larger than 600 EDs. This is because as the instance size grows, more ED itineraries with overlap are generated. Then the EP fleet gets more local or distant switch opportunities, and has more options in selecting charging customers. In other words, a more flexible charging schedule renders each EP a higher service rate. Those observations provide valuable insights. We can foresee that the EV usage and the service efficiency of the CaaS platform present positive interaction and encourage each other's growth in a certain range of EV market penetration.

Even though the results in Table 1 shows the growth of the EV market (i.e., large ED instances) will encourage a better service efficiency of the CaaS platform (i.e., a large service rate per EP), it is also noticed that the EP fleet size nonlinearly increases with the ED demand size, which is further subject to the total traffic demand (traffic condition) and the EV penetration. As the usage of EV grows, we can predict that a large EP fleet size is needed, which may inevitably worsen traffic conditions. Therefore, it is critical to investigate the traffic impact coming from the CaaS platform under different EV market penetrations and suggest the application environment. Accordingly, this subsection aims to answer the question: under a certain traffic condition, what is the range of the EV market penetration for the CaaS platform to operate acceptable, while its traffic overhead is of critical concern?

To do that, our experiments consider three traffic scenarios defined by the network volume to capacity: v/c ratio equal to 0.35 (low or no congestion), 0.6 (moderate congestion), and 0.8 (heavy congestion) (Xie et al., 2018) respectively. They can be generated by randomly selecting the background traffic from the realistic Chicago O-D database until the desired v/c ratio is reached. For each scenario, we simulate 10 instances for different EV penetration rates. All parameters follow the default setting. Specifically, each EP is equipped with two 90 kWh batteries, e.g., the Tesla P90D battery. As 540 kg battery weight gives a total weight of 1080 kg, the EP is a small van with one passenger car unit (PCE) (Adnan, 2014).

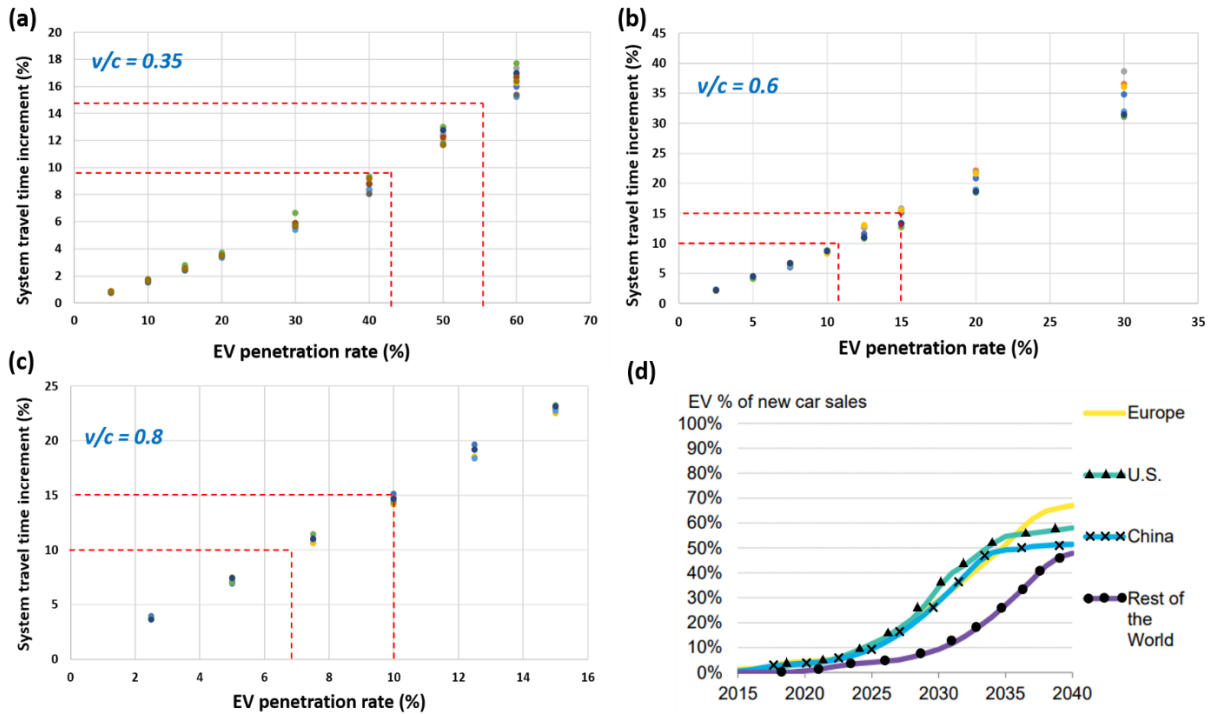

**Figure 8.** Traffic overhead vs. EV penetration under three traffic scenarios: (a) low or no congestion; (b) moderate congestion; (c) heavy congestion. And (d) long-term EV sales penetration prediction

(Source: Bloomberg New Energy Finance)

The traffic impact of the CaaS platform is quantified by comparing the system travel time with or without involving the EP fleet. Briefly, a proportion of background traffic is selected as the EDs according to the EV market penetration rate. Then, a fleet of EPs are dispatched according to the charging plans from the CCWS algorithm. After that, we can update the link volume and compute the system travel time for the vehicles (except EPs) as $STT_E$. Then, the traffic impact of the CaaS platform is assessed by the increase of system travel time measured by Equation (40).

$$\Delta STT = \frac{STT_E - STT_0}{STT_E} \times 100\%, \tag{40}$$

where $STT_0$ and $STT_E$ respectively represent the system travel time without/with dispatching EP fleets. The link travel time is measured by the BPR function with corresponding link volume.

Figure 8 (a)-(c) display $\Delta STT$ under different EV penetration rate and corresponding maximum EV penetration under a tolerable $\Delta STT$ such as 10 % or 15%. It means that a delay less than 2 minutes (for 10 % $\Delta STT$) or 3 minutes (for 15 % $\Delta STT$) is not significant for a 20 minutes trip.

Accordingly, Figure 8 (a) indicates that the traffic overhead resulting from the CaaS platform is acceptable – $\Delta STT$ is less than 15% (or 10%) – when traffic is under light traffic congestion (e.g., 0.35 v/c) and EV penetration is below 55% (or 42%). As the network traffic becomes more congested, such as under 0.6 v/c in Figure 8 (b), the CaaS platform is only suitable to EV penetration market below 15% (or 11%) for the tolerable $\Delta STT$ equal to 15% (or 10%). Figure 8 (c) demonstrates the consistent trend as the initial traffic environment becomes more congested. In conclusion, the results in Figure 8 (a)-(c) indicate that the CaaS platform operates well and won't introduce severe traffic overhead under light traffic condition with low EV penetration (less than 50%), but its usage should be limited under congested traffic condition such as in CBD area or during peak hours. Moreover, Figure 8(d) shows that EV penetration (i.e., EV%) is less than 60% in most of the countries in the year 2040 (Finance, 2017). Therefore, we can foresee that the CaaS platform will be an applicable solution for the next few decades to cultivate the mass-adoption of EVs.

### 5.1.4 EP energy inventory analysis

In practice, each EP can carry multiple batteries. Thus, individual EPs may have different initial energy inventories. They together will affect the service efficiency of the CaaS platform. It is valuable to test the sensitivity of the platform performance to the initial energy inventory. The obtained insights can help the CaaS platform develop strategic energy inventory control and improve the service efficiency. To set up the experiment, we randomly generate and simulate 10 instances, each with an instance size of 10,000 EDs. We collect the EP service rate under the experiments with different initial energy inventories, varying from 36 kWh to 243 kWh, and then plot the results in Figure 9. Note that the capacity of one battery is 90 kWh.

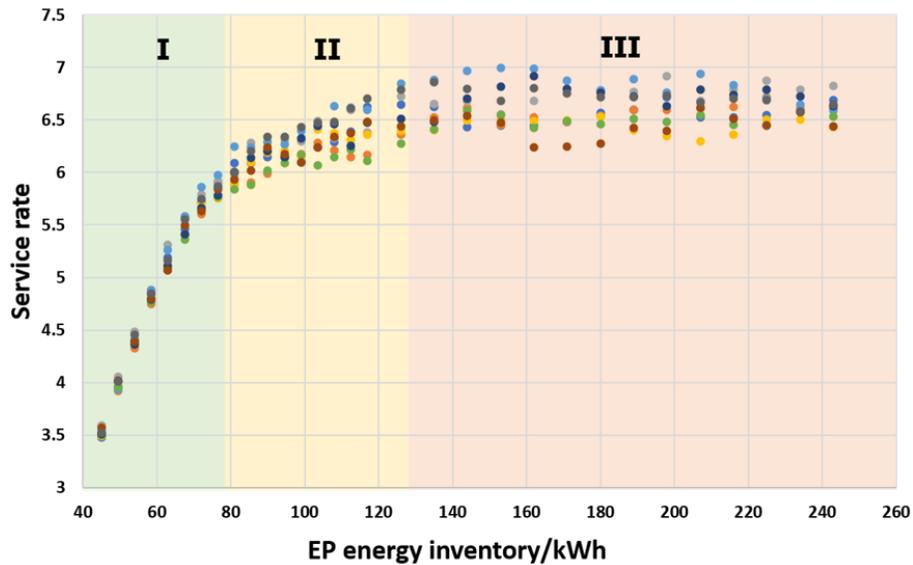

**Figure 9.** Plot of EP service rate under the instances with different initial EP energy inventory.

The results in Figure 9 demonstrate that the service rate first grows with the increase of individual EP's initial energy inventory and then reaches a plateau when the inventory is above 130 kWh (about 1.4 battery). More exactly, the results present three phases (I-III). When the initial battery is less than 80 kWh, the service rate almost linearly grows with the increase of EP's initial energy inventory; after that, the effect slows down, and eventually reaches a plateau when we further increase the initial battery inventory above 130kWh. This implies that the initial energy inventory only limits the service rate of an EP when it is relatively low. There are other factors to affect the service rate. We noticed that the EDs' itineraries and their inherent correlations play another important factor. For example, each ED service sets a given itinerary and service time window. Due to this limit, an EP may not be able to reach an ED timely even if it has enough energy inventory. Therefore, there exists an upper bound for the EPs' service rate as it presents in the phase III. In the phase II, the service rate is limited by mixed factors and therefore, the plot exhibits a slower increase. The factor of the EDs' itineraries and their inherent correlations can also be evidenced by the much more spreading data points in the phase II and III.

More importantly, the results in Figure 9 imply that an EP is ready to go if its initial energy inventory is above 130 kWh. In other words, an EP does not have to be fully recharged for the next duty. Consequently, an EP can participate the service in a new batch from the location in the previous service cycle without returning to a depot, if its remaining inventory is above 130 kWh. Accordingly, we can foresee that if each EP carries more than 2 batteries, it can reduce the times for returning to the depot and improve the service coverage by obtaining more opportunities to timely serve those EDs that are far away from depots, but close to the idle EPs from the last service cycle. This insight can help set up a more efficient CaaS platform and we have factored in these considerations in our model by allowing the EPs to start from the current location (not limited to the depots) with the current energy inventory.

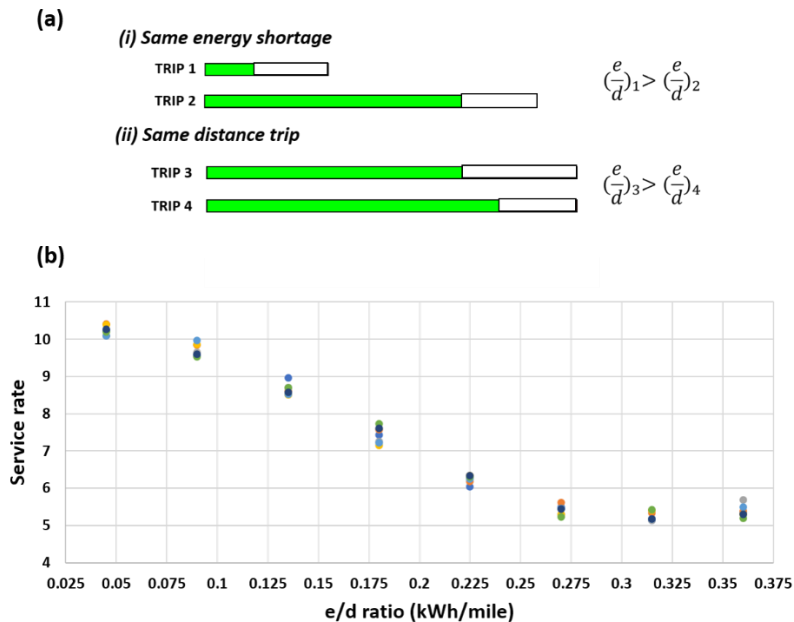

**Figure 10.** (a) Schematic representation of the energy shortage of EDs; the length of each bar represents the length of the trip. The filled/unfilled section represents the remain/shortage of battery; (b) EP service rates under different ED energy shortage rate.

### 5.1.5 Pricing strategy

This study noticed that the energy shortages and trip lengths of the EDs jointly affect the service efficiency of the CaaS platform. We demonstrate this observation by using the trips in Figure 10 (a) as examples. Given Trips 1 and 2 have the same energy shortage but different trip lengths. Our experiments showed that Trip 2 will potentially give more flexibility to schedule an EP trip than Trip 1 does since Trip 2 is longer and it can accommodate more chances to encounter an EP by either local or distant switch. On the other hand, Trips 3 and 4 have the same trip length, but Trip 4 needs less energy. Our experiments show that it is relatively easier to consolidate the service for Trip 4 with others served by an EP. Thus, Trip 4 will facilitate the service schedule better than Trip 3. With this view, we define the energy shortage rate (labeled as e/d ratio) by the ratio of the ED's energy shortage to its trip length, and then use the e/d ratio to quantify the collective effect. On the other hand, this e/d ratio provides a quantitative way to measure the range anxiety. Namely, a larger value of the e/d ratio indicates a higher level of range anxiety.

Accordingly, we conduct experiments to explore comprehensive insights for the effect of the e/d ratio on the service efficiency of the CaaS, in particular, the service rate. To do that, our experiments randomly generate 10 instances, each of which involves 10,000 EDs with an average e/d ratio such as 0.05 kWh/mile. The 10 instances cover the scenarios that the e/d ratio varies from 0.025 kWh/mile to 0.375 kWh/mile. All other parameters follow the default setting. From the results in Figure 10 (b), we noticed that the EP service rate gets worse as the e/d ratio increases, but becomes stable when the e/d ratio is above 0.25 kWh/mile. This interesting observation gives a hint for the CaaS platform to improve the service rate by developing a pricing strategy towards the ED's e/d ratio. For example, the CaaS may charge more if an ED is with a higher e/d ratio. Then the proposed pricing strategy encourages the EDs to request the mE2 charging service when they have a small e/d ratio, and ultimately mitigate the range anxiety.

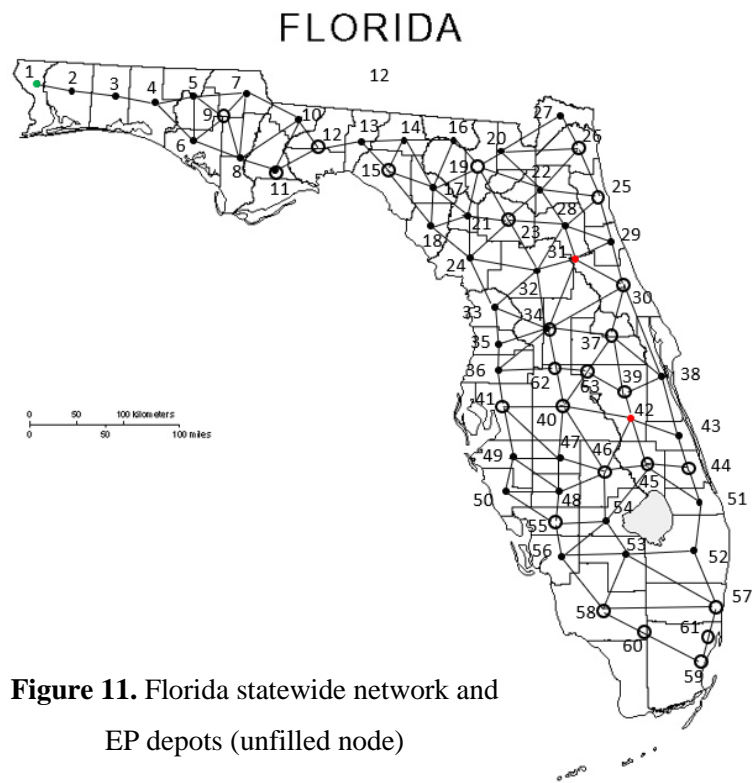

**Figure 11.** Florida statewide network and EP depots (unfilled node)

## 5.2 Numerical study 2: Florida statewide network

Numerical study 1 evaluates the performance of the CaaS platform on a citywide network. The ED trips are all intra-city trips with 11 miles as the average trip length. However, existing studies show that the statewide long-distance trips usually suffer from the range anxiety more due to the battery range limitation and sparse charging service stations on the inter-city roads. With this view, the CaaS service can be a promising relief for long-distance travelers using EVs. This study therefore further tests the performance of the CaaS platform at a state-level network.

We use Florida statewide network as the testbed. As shown in Figure 11, the county-based network is with 63 nodes representing counties, big cities, important road junctions, or EP depots, and 260 links representing the connections. The potential EP recharging depots (circled nodes) are selected based on the existing DC fast charging (DCFC) locations (FDOT, 2020). The statewide freight trips were obtained by running the simulation model, Florida Statewide Model v7.2 (FLSWM) (Zanjani et al., 2015). The average passenger trip distance is 231.3 miles, with the minimum of 22.5 miles and the maximum of 937.6 miles.

To analyze the impact of the trip distance on the performance of the CaaS platform, we categorize the passenger trips into seven classes in terms of their trip distances, $[100,150)$, $[150,200)$, $[200,250)$, $[250,300)$, $[300,350)$, $[350,400)$, $[400,450)$. For each distance class, we generate 10 instances by selecting 10,000 qualified trips from the passenger trip dataset. All parameters follow the default setting, except that the EPs can carry either 2 batteries (180 kWh capacity) or 6 batteries (540 kWh capacity).

Figure 12 displays the plots of EP service rate vs. the ED trip distance under different two scenarios, in which an EP is equipped with 2 or 6 batteries, respectively. The results indicate that the max EP service rate is around 2.5 and it declines further as the length of individual ED's trip increases. Equipped with more batteries helps to improve the service efficiency. However, a six-battery pack gives a total weight of 3,240 kg, so the EP has to be a light-duty truck with 3 PCE (Adnan, 2014). Thus, we conclude that using the CaaS platform to serve state-wide trips may cause more traffic than it does for citywide trips. This is consistent with our common sense. The long-distance EDs (such as trucks) usually require more energy, while each EP can only carry limited batteries. Moreover, the EPs usually travel longer distances to provide statewide service and then return the depots. Therefore, they experience more idle energy loss. In conclusion, these results imply that building more charging stations (depots) will benefit the CaaS platform on the statewide service. Moreover, the results also suggest using a pricing strategy according to the EDs' trip lengths and energy requests will balance the service quality, operation cost, and system traffic impacts.

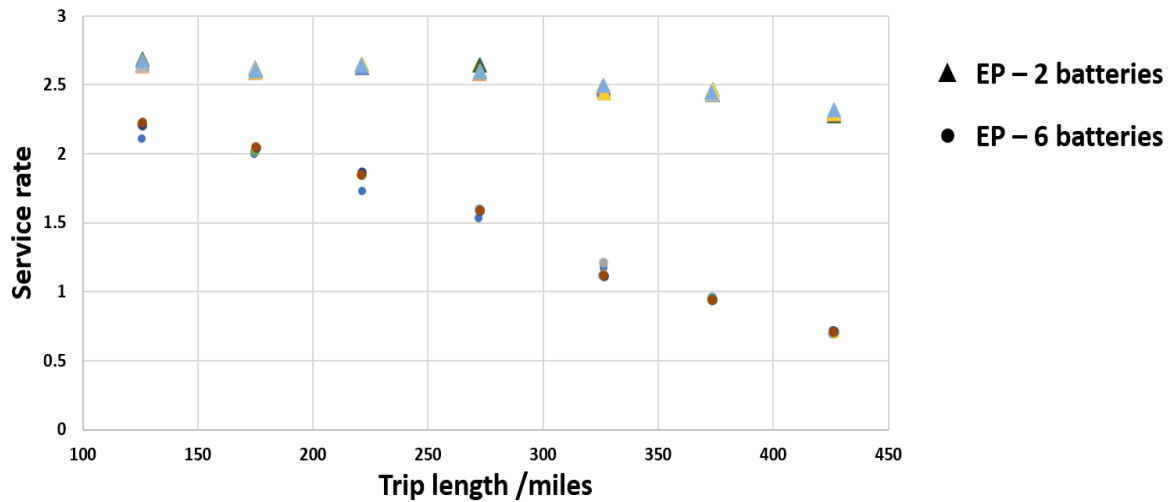

**Figure 12.** Plot of EP service rate under the instances with different ED trip distance and each EP carrying different numbers of batteries.

## 6. Conclusion

This paper develops a CaaS platform that provides the charging services for EVs on the move so that we can relieve the range anxiety and charging delay of EDs by taking advantage of the mE2 technology. Mathematically, the CaaS platform is modeled as a new mE2-VRP model, which is a mixed-integer program that seeks to explore the optimal routing schemes so that we can dispatch the commercial electricity suppliers (i.e., EPs) to charge EVs with the low level of batteries (i.e., EDs) when both of them are on the move while minimizing the extra traffic introduced by EPs. The model development addressed the challenges raised by the unique features of the on-the-move electricity delivery. Moreover, we contribute a new heuristic algorithm (i.e., the CCWS algorithm) to efficiently solve a large-scale mE2-VRP by integrating clustering and parallel computation with commercial optimization solver. The efficiency and applicability of our approaches are accessed by the numerical experiments built upon the Chicago city network and Florida state network. The main results and findings are summarized as follows.

First of all, we examined the performance of the CCWS algorithm. The experiments built upon both testbeds show that the CCWS algorithm outperforms commercial solvers such as Gurobi. It addresses the scalability of the mE2-VRP model and enables us to evaluate the performance of the CaaS platform with realistic settings under different traffic circumstances. We next accessed the performance of the CaaS in the Chicago city network. The experiment results indicate that the CaaS platform introduces extra traffic by dispatching the EPs. This traffic overhead is acceptable under the light traffic environment with low EV market penetration (less than 50%). Thus, the CaaS platform is a promising and applicable solution to cultivate the mainstream adoption of EVs in the near future. Last, we conducted sensitivity analyses on the operational inputs of the CaaS platform. The results demonstrate that the CaaS platform may maintain a high service efficiency by setting the initial energy inventory of the EPs above 130 kWh since it avoids frequent recharging trips. More importantly, the CaaS platform may lead to less traffic overhead if the majority of the EDs are served when they are under a lower energy shortage. The CaaS platform performs better for the EDs with shorter intercity trips than those conducting longer statewide trips. These observations together suggest using pricing schemes according to the EDs' energy requests and trip lengths to balance the service quality, operation cost, and system traffic impacts.

This study is among the early efforts to study the on-the-move charging service to relieve range anxiety and promote the usage of EVs. Several potential future research can be further explored from this study. For example, we will investigate the operation of the CaaS platform while integrates the depot location planning and pricing strategies. We noticed that travel time and energy usage often present stochastic features. Our future work will consider these factors and make the mE2-VRP can accommodate the uncertainty in reality. Those extensions will make the CaaS work more efficiently, but also complicate the existing models and raise new computation challenges. Our future work will address them. This CaaS platform promises online services, but the mE2-VRP model is a large-scale mixed-integer program, which is hard to solve efficiently in general. To address this difficulty, the CCWS algorithm adapts the CaaS platform to provide timely online services under a large-scale instance in practice.

**Acknowledgment**
This research is partially supported by National Science Foundation award CMMI 1818526.

# Appendix
## Appendix A. Nomenclature

**Table 2:** Parameters and Variables

| 1. Parameters in the mE2-VRP | |
|---|---|
| Notation | Explanation |
| $H = \{h\}_{h=1}^{H}$ | One planning horizon, discretized into $H$ time intervals. |
| $T > 0$ | The length of the time interval. |
| $i \in \mathcal{D}_h$ | Index and set (abbreviated as $\mathcal{D}$) of the EDs that depart during the $h$ time interval. |
| $s \in \mathcal{S}_h$ | Index and set (abbreviated as $\mathcal{S}$) of the EPs available at beginning of the $h$ time interval. |
| $D_n$ | A set of the EDs that pass through a node $n$. |
| $\mathcal{N}$ | The set of nodes. $\mathcal{N} = N_0 \cup N_e \cup N_p \cup N$. |
| $N_0$ | A set of the first available locations for service along the routes of the EDs. |
| $N_e$ | A set of the destinations of the EDs. |
| $N_p$ | A set of depots of the EPs. |
| $n \in N$ | Index and set of encounter nodes between the EDs. |
| $N_i$ | A set of encounter nodes along the route of ED(i). |
| $\mathcal{N}_i$ | A set of nodes to represent the route of ED(i), $\mathcal{N}_i = (i_0, N_i, i_e)$. |
| $A_i \subseteq A$ | A set of arcs to represent the route of ED(i). |
| $a_n^i$ | The arc starting from node $n$ on the route of ED(i). $a_n^i$ is replaced by $a^i(n)$ in some formulations. |
| $i_0 \in N_0$ | Index and set of origin of the EDs. |
| $i_e \in N_e$ | Index and set of destination of the EDs. |
| $s_0$ | Index of depots/current locations of the idle EPs. |
| $p_n$ | Index of the nearest depot from node $n$. $p_n \in N_p$. |
| $\tilde{t}_0^i$ | The estimated arrival time of ED(i) at $i_0$. |
| $\bar{\tau}_i$ | The maximum endurable waiting time of ED(i) at $i_0$. |
| $\bar{e}_i$ | The battery capacity of ED(i). |
| $e_0^i$ | The estimated battery level of ED(i) at $i_0$. |
| $\varpi_i$ | The energy consumption rate of ED(i). |
| $\underline{e}_d$ | The safety energy inventory of the EDs. |
| $e_{i,a}^-$ | The energy that ED(i) losses on the arc $a$. |
| $\bar{e}_{i,a}^+$ | The maximum energy that ED(i) can receive from the EPs on arc $a$. |
| $e_0^s$ | The initial energy inventory of EP(s). |
| $\varpi_s$ | The energy consumption rate of ED(s). |
| $\underline{e}_s$ | The safety energy inventory of the EPs. |
| $e_{n \to m}^-$ | The energy consumed for the EPs to travel from node $n$ to $m$. |
| $\eta$ | The mE2 power transfer rate of the EPs. |
| $e_{s,a}^-$ | The energy that EP(s) losses on the arc $a$. |
| $t_{n \to m}$ | The average travel time from node $n$ to node $m$. |
| **2. Variables in the mE2-VRP** | |
| $F$ | Optimization objective of the mE2-VRP. |
| $z_{i,a}^s$ | Binary variable to identify whether EP(s) serves ED(i) on the arc $a$. If true, $z_{i,a}^s = 1$, otherwise $z_{i,a}^s = 0$, $\forall s \in \mathcal{S}, i \in \mathcal{D}, a \in A_i$. $Z = \{z_{i,a}^s\}$ is the set. |
| $o_i^{n,s}$ | Binary variable to identify whether EP(s) travels from depot to charge the ED(i) at node $n$. If true, $o_i^{n,s} = 1$, otherwise $o_i^{n,s} = 0$, $\forall s \in \mathcal{S}, i \in \mathcal{D}, n \in \mathcal{N}_i \backslash i_e$. $O = \{o_i^{n,s}\}$ is the set. |
| $q_i^{n,s}$ | Binary variable to identify whether EP(s) leaves ED(i) and return to depot from node $n$. If true, $q_i^{n,s} = 1$, otherwise $q_i^{n,s} = 0$, $\forall s \in \mathcal{S}, i \in \mathcal{D}, n \in \mathcal{N}_i \backslash i_0$. $Q = \{q_i^{n,s}\}$ is the set. |
| $u_{i,j}^{n,s}$ | Binary variable to identify whether EP(s) conducts local switches from ED(i) to ED(j) at node $n$. If true, $u_{i,j}^{n,s} = 1$, otherwise $u_{i,j}^{n,s} = 0$, $\forall s \in \mathcal{S}, i \in \mathcal{D}, n \in \mathcal{N}_i \backslash i_0, j \in D_n, j \neq i$. $U = \{u_{i,j}^{n,s}\}$ is the set. |
| $v_{i,j}^{n,m,s}$ | Binary variable to identify whether EP(s) conducts distant switch from ED(i) at node $n$ to ED(j) at node $m$. If true, $v_{i,j}^{n,m,s} = 1$, otherwise $v_{i,j}^{n,m,s} = 0$, $\forall s \in \mathcal{S}, i, j \in \mathcal{D}, i \neq j, n \in \mathcal{N}_i \backslash i_0, m \in \mathcal{N}_j \backslash j_e, n \neq m$. $V = \{v_{i,j}^{n,m,s}\}$ is the set. |
| $w_i^{n,s}$ | Binary variable to identify whether EP(s) charges ED(i) on the arc $a(i,n)-1$ and continue to charge ED(i) on the arc $a^i(n)$. If true, $w_i^{n,s} = 1$, otherwise $w_i^{n,s} = 0$, $\forall s \in \mathcal{S}, i \in \mathcal{D}, n \in N_i$. $W = \{w_i^{n,s}\}$ is the set. |
| $e_{i,a}^+$ | Continuous variable represents the energy that ED(i) receives from the EPs on arc $a$, $\forall i \in \mathcal{D}, a \in A_i$. |
| $e_n^i$ | Continuous variable represents the energy inventory of ED(i) at node $n$, $\forall i \in \mathcal{D}, n \in \mathcal{N}_i \setminus i_0$. |
| $\tau_i$ | Continuous variable represents the waiting time of ED(i) at $i_0$. $\forall i \in \mathcal{D}$. |
| $t_0^i$ | Continuous variable represents the actual departure time of ED(i), $\forall i \in \mathcal{D}$. $t_0^i = \tilde{t}_0^i + \tau_i$. |
| $t_n^i$ | Continuous variable represents the arrival time of ED(i) at node $n$, $\forall i \in \mathcal{D}, n \in \mathcal{N}_i$. $t_n^i = \tilde{t}_0^i + t_{i_0,n} + \tau_i$. |
| $E_s$ | Continuous variable represents the total energy consumed for EP(s), $\forall s \in \mathcal{S}$. |
| $\gamma_{i,a}$ | Continuous auxiliary variable to linearize the non-linear term of $e_{i,a}^+(\sum_{s \in \mathcal{S}} z_{i,a}^s)$, $\forall i \in \mathcal{D}, a \in A_i$. |
| $\zeta_{i,a}$ | Continuous auxiliary variable to linearize the non-linear term of $z_{i,a}^s e_{i,a}^+$, $\forall i \in \mathcal{D}, a \in A_i$. |
| **3. Parameters and variables in the CCWS algorithm** | |
| $J_{i,n}$ | Local-Switch Candidate set of ED(i) at node $n$ $\forall i \in \mathcal{D}, n \in \mathcal{N}_i \setminus i_0$. |
| $p(i,j)$ | Number of times that an EP can do local switch from ED(i) to ED(j) along the route of ED(i). |
| $P$ | Opportunity matrix, $P = [p(i,j)]_{i,j \in D}$. |
| $\mathcal{S}_0$ | A set of EPs that are dispatched in the initial solution. |
| $\mathcal{D}_s^*$ | The set of the EDs serviced by EP(s) in the initial feasible solution. |
| $s(n)$ | The $n$-th ED served by EP(s) in the initial feasible solution. $s(-1)$ is the last ED served by EP(s). |
| $E_s^*$ | The total energy loss (exclude idle energy loss from depot to the first customer and from the last customer to the depot) of EP(s) in the initial feasible solution. |
| $t_0^{*s}$ | The time that EP(s) starts charging the first customer in the initial feasible solution. |
| $t_e^{*s}$ | The time that EP(s) leaves the last customer in the initial feasible solution. |
| $H(s)$ | The seed tour indexed by the EP(s) in the initial feasible solution, $\forall s \in \mathcal{S}_0$, which is a tuple, $H(s) = (n, t_0^{*s}, m, t_e^{*s}, E_s^*)$, and records the locations and time that EP(s) starts servicing the first customer and completing service with the last customer as well as the total energy loss. |

## Appendix B. Algorithm 1: Clustering Algorithm

**Algorithm 1** Clustering Algorithm

---

Initialize: cluster index $c \leftarrow 0$
1   **While** ( $p(i,j) \neq 0, \forall i,j$ ) **do**
2       $c \leftarrow c + 1$
2       $i \leftarrow argmax_i[p(i,j)]$, Construct cluster $\mathcal{C}(c) = \{i\}$; $p(j,i) = 0, \forall j$.
3       $j \leftarrow argmax_j[p(i,j)]$
4       **While** ( $p(i,j) > 1$ ) **do**
5           $\mathcal{C}(c) \leftarrow \mathcal{C}(c) \cup \{j\}$
6           $p(k,j) = 0, \forall k$
7           $i \leftarrow j$; $j \leftarrow argmax_k[p(i,k)]$
8   **return** $\mathcal{C}$

---

The algorithm starts with the largest entry $p(i,j)$ in the opportunity matrix P to construct the first cluster $\mathcal{C}(1) = \{i,j\}$. ED(i) and ED(j) are in one cluster because ED(j) has a great potential to be served next by an EP doing local switch from ED(i). Now as both ED(i) and ED(j) are clustered, they cannot present in other clusters. Accordingly, the algorithm zeros the $i$-th and $j$-th column of the opportunity matrix. Next, the algorithm moves to the row $j$, and searches for the largest entry to include the next ED to the cluster $\mathcal{C}(1)$. As the largest entry is searched, this ED has the highest potential to be serviced after ED(j). As the algorithm repeats the searching process, the cluster $\mathcal{C}(1)$ is settled until no more ED can be clustered into $\mathcal{C}(1)$. Then, the algorithm searches the largest entry in the updated opportunity matrix and repeats the whole procedure to construct other clusters. Finally, the clustering algorithm returns a set of clusters, $\mathcal{C}(c), \forall c$. For convenience, we denote the set of EDs that are associated with the cluster $\mathcal{C}(c)$ as $\mathcal{D}_c$.

## Appendix C. Algorithm 2: Merging Algorithm

**Algorithm 2** Merging Algorithm

| | |
|---|---|
| 1 | Create initial seed tours, $H(\mathcal{S}_0) = (H(s), \forall s \in \mathcal{S}_0)$ by solving sub-mE2-VRP$(c), \forall c$. |
| 2 | **for** $s \in \mathcal{S}_0$ **do** |
| 3 |   **for** $l \in \mathcal{S}_0$ **do** |
| 4 |     **if** $H(s) \vee H(l)$ **then*** |
| 5 |       $l^* \leftarrow l;\ t_e^{s-l} \leftarrow \max\{t_e^{*s} + t_{n_2, m_1(l)}, t_0^{*l}\} + t_e^{*l} - t_0^{*l}$ |
| 6 |       **for** $k \in \mathcal{S}_0$ **do** |
| 7 |         **if** $H(s) \vee H(k)$ **and** $t^* > \max\{t_e^{*s} + t_{n_2, m_1(k)}, t_0^{*k}\} + t_e^{*k} - t_0^{*k}$ **then** |
| 8 |           $l^* \leftarrow k;\ t_e^{s-l^*} \leftarrow \max\{t_e^{*s} + t_{n_2, m_1(k)}, t_0^{*k}\} + t_e^{*k} - t_0^{*k}$ |
| 9 |       Merge $H(s)$ with $H(l^*)$; $H(s) \leftarrow \bar{H}(s - l^*)$ |
| 10 |       $\mathcal{S}_0 \leftarrow \mathcal{S}_0 \setminus \{l^*\}$ |
| 11 | **return** $H(\mathcal{S}_0)$ |

\* $H(s) \vee H(l)$ indicates that tour $H(s)$ can merge with tour $H(l)$ by EP(s).

Mainly, a set of seed tours $H(\mathcal{S}_0) = (H(s), \forall s \in \mathcal{S}_0)$ is solved from sub-mE2-VRP$(c)$ (Line 2). Next, for any tour $H(s) \in H(\mathcal{S}_0)$, we consider tour $H(l^*)$ has the highest priority to be merged with $H(s)$ by EP(s), if $H(l^*)$ leads to the smallest departure time in the merged tour $\bar{H}(s - l^*)$ (Line 4 – 9). For example, the time that EP separates from the last customer is $t_e^{s-l} = \max\{t_e^{*s} + t_{n_2, m_1(l)}, t_0^{*l}\} + t_e^{*l} - t_0^{*l}$ in the merged tour $\bar{H}(s - l) = H(s) \vee H(l)$. The algorithm compares $t_e^{s-l}$ when $H(s)$ merges with different $H(l) \in H(\mathcal{S}_0)$ and selects $H(l^*)$ with the smallest $t_e^{s-l^*}$ to be merged with $H(s)$. After merging two seed tours, we update the tour information as $\bar{H}(s - l^*) = (n_1, t_0^{*s}, m_2, t_e^{s-l^*}, \bar{E}_s^*)$, where $\bar{E}_s^* = E_s^* + e_{n_2 \to m_1}^- + E_l^*$, which is the total energy consumed for EP(s) to serve all demands in both seed tours. And the departure time, $t_e^{s-l^*}$ is updated by the Equation (41).

$$t_e^{s-l^*} = \max\{t_e^{*s} + t_{n_2 \to m_1}, t_0^{*l}\} + t_e^{*l} - t_0^{*l} \tag{41}$$

Also,

$$\Delta \tau_i = \begin{cases} 0, & \text{if } t_e^{*s} + t_{n_2 \to m_1} \leq t_0^{*l} \\ t_e^{*s} + t_{n_2 \to m_1} - t_0^{*l}, & \text{if } t_e^{*s} + t_{n_2 \to m_1} > t_0^{*l} \end{cases} \tag{42}$$

Finally, the tour $H(l^*)$ is removed from tours set $H(\mathcal{S}_0)$ and the process is repeated until no more tours can be merged.

## Appendix D. Proof of Lemma 2 and 3
**Lemma 2.** *For any ED(i), a trip satisfying (13) and (14) eliminates subtours that are only formed by local switch from/back to ED(i).*

*Proof.* We also prove Lemma 2 by showing contradiction. Similar with the proof of Lemma 1, we first show that constraints (13) and (14) eliminate an illegal subtour formed by multiple local switches, when it only involves 2 EDs ($n = 2$). Then, the proof is extended to the case involving more than 2 EDs ($n \geq 3$). Suppose there exists an illegal subtour purely formed by local switches with n = 2. Figure 13(a) shows an example. An EP first serves $ED(i_1)$ as they move together from node $n_1$ to node $n_2$, and then it locally switches to serve $ED(i_2)$ starting from node $n_2$ (local switch: $u_{i_1,i_2}^{n_2,s} = 1$) until arriving at node $n_1$. After that, the EP revisits $ED(i_1)$ at node $n_1$ (local switch: $u_{i_2,i_1}^{n_1,s} = 1$). We have (43)(44) below, given that constraint (13) is satisfied with $u_{i_1,i_2}^{n_2,s} = u_{i_2,i_1}^{n_1,s} = 1$,

$$t_{i_2}^{n_2} - t_{i_1}^{n_2} \geq 0 \tag{43}$$

$$t_{i_1}^{n_1} - t_{i_2}^{n_1} \geq 0 \tag{44}$$

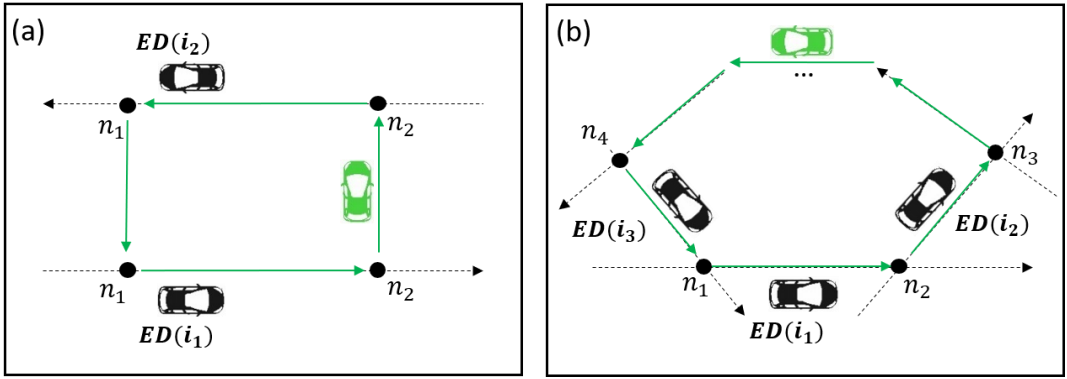

**Figure 13.** illegal subtour with size (a) $n = 2$ and (b) $n \geq 3$, that is formed by local switches only.

Adding (43)(44), we have (45) below.

$$0 > \left(t_{i_2}^{n_2} - t_{i_2}^{n_1}\right) + \left(t_{i_1}^{n_1} - t_{i_1}^{n_2}\right) \geq 0 \tag{45}$$

which leads to the contradiction. Therefore, such subtour does not exist in a feasible solution of mE2-VRP because constraint (13) will be violated.

When $n \geq 3$ (see Figure 13 (b)), we assume an EP serves $ED(i_1)$ from node $n_1$ to node $n_2$, and then it locally switches to serve $ED(i_2)$ starting from node $n_2$ (local switch: $u_{i_1,i_2}^{n_2,s} = 1$) to node $n_3$. Next, it locally switches to serve other EDs until it serves $ED(i_3)$ from node $n_4$ to node $n_1$. And finally, the EP revisits $ED(i_1)$ at node $n_1$ (local switch: $u_{i_3,i_1}^{n_1,s} = 1$). To show such subtour does not exist in a feasible solution, we only need to prove at least one of the local switches involved in the subtour is infeasible (violates constraint (13)). Without loss of generality, we only focus on two local switches from $ED(i_1)$ to $ED(i_2)$ at node $n_2$, and $ED(i_3)$ to $ED(i_1)$ at node $n_1$ and assume all other local switches are feasible. We have (46)(47) below, given that constraint (13) is satisfied with $u_{i_1,i_2}^{n_2,s} = u_{i_3,i_1}^{n_1,s} = 1$,

$$t_{i_2}^{n_2} - t_{i_1}^{n_2} \geq 0 \qquad (46)$$

$$t_{i_1}^{n_1} - t_{i_3}^{n_1} \geq 0 \qquad (47)$$

Adding (46)(47), we have (48) below,

$$0 \geq \left(t_{i_1}^{n_2} - t_{i_1}^{n_1}\right) + \left(t_{i_3}^{n_1} - t_{i_2}^{n_2}\right) > 0 \qquad (48)$$

which leads to a contradiction. The subtour does not exist, which completes the proof.

**Lemma 3.** *For any ED(i), a trip satisfying (13) and (14) eliminates subtours that are formed by distant switch from/back to and local switch back to/from ED(i).*

*Proof.* We also prove Lemma 3 by showing contradiction. Similar with the proof of Lemma 1 and 2, we first show that constraints (13) and (14) eliminate an illegal subtour formed by multiple local and distant switches, when it only involves 2 EDs ($n = 2$). Then, the proof is extended to the case involving more than 2 EDs ($n \geq 3$). Suppose there exists an illegal subtour formed by distant and local switches with $n = 2$, see Figure 14(a). Both ED($i_1$) and ED($i_2$) travel from node $n_1$ to node $n_2$. An EP first serves ED($i_1$) from node $n_1$ to node $n_2$, and then it travels back to node $n_1$ and conducts distant switch ($v_{i_1,i_2}^{n_2,n_1,s} = 1$) to serve ED($i_2$) from node node $n_1$ until arriving at node $n_2$. After that, the EP revisits ED($i_1$) at node $n_2$ (local switch: $u_{i_2,i_1}^{n_2,s} = 1$). We have (49)(50) below, given that constraint (13) and (14) are satisfied with $u_{i_2,i_1}^{n_2,s} = v_{i_1,i_2}^{n_2,n_1,s} = 1$,

$$t_{i_2}^{n_1} - t_{i_1}^{n_2} - t_{n_2 \rightarrow n_1} \geq 0 \qquad (49)$$

$$t_{i_1}^{n_2} - t_{i_2}^{n_2} \geq 0 \qquad (50)$$

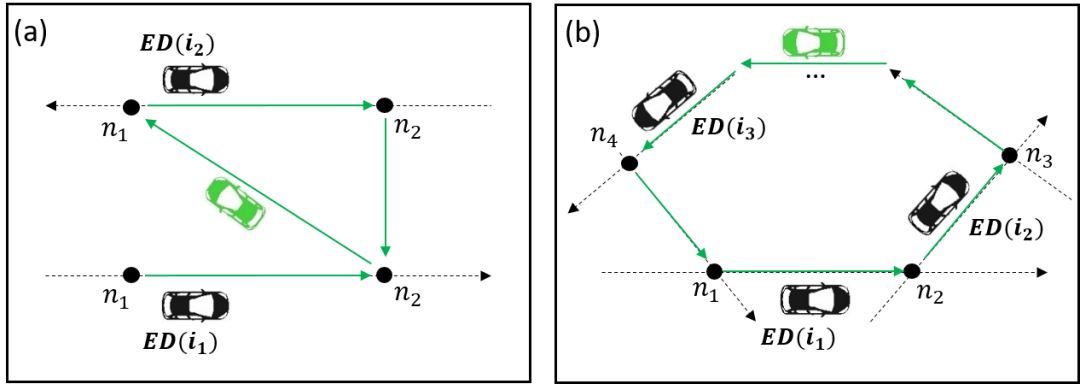

**Figure 14.** illegal subtour formed by distant and local switches with size (a) $n = 2$ and (b) $n \geq 3$.

Adding (49)(50), we have (51) below,

$$-t_{n_2 \to n_1} \geq \left(t_{i_1}^{n_2} - t_{i_2}^{n_1}\right) + \left(t_{i_2}^{n_2} - t_{i_2}^{n_1}\right) > 0 \tag{51}$$

(51) indicates that $t_{n_2 \to n_1} < 0$ which contradicts with $t_{n_2 \to n_1} > 0$. Therefore, such subtour does not exist in a feasible solution of mE2-VRP because constraints (13) and (14) will be violated.

When $n \geq 3$ (see Figure 14(b)), we assume an EP serves $ED(i_1)$ from node $n_1$ to node $n_2$, and then it locally switches to serve $ED(i_2)$ starting from node $n_2$ (local switch: $u_{i_1,i_2}^{n_2,S} = 1$) to node $n_3$. Next, it conducts local or distant switches to serve other EDs until it serves $ED(i_3)$ and arrives at $n_4$. And finally, the EP revisits $ED(i_1)$ at node $n_1$ by distant switch from node $n_4$ to node $n_1$ (distant switch: $v_{i_3,i_1}^{n_4,n_1,S} = 1$). To show such subtour does not exist in a feasible solution, we only need to prove at least one of the local or distant switches involved in the subtour is infeasible (violates constraints (13) and (14)). Without loss of generality, we only focus on the local switch from $ED(i_1)$ to $ED(i_2)$ at node $n_2$, and distant switch from $ED(i_3)$ at node $n_4$ to $ED(i_1)$ at node $n_1$ and assume all other switches are feasible. We have (52)(53) below, given that constraint (13) and (14) are satisfied with $u_{i_1,i_2}^{n_2,S} = v_{i_3,i_1}^{n_4,n_1,S} = 1$,

$$t_{i_2}^{n_2} - t_{i_1}^{n_2} \geq 0 \tag{52}$$

$$t_{i_1}^{n_1} - t_{i_3}^{n_4} - t_{n_4 \to n_1} \geq 0 \tag{53}$$

Adding (52)(53), we have (54) below,

$$-t_{n_4 \to n_1} \geq \left(t_{i_1}^{n_2} - t_{i_1}^{n_1}\right) + \left(t_{i_3}^{n_4} - t_{i_2}^{n_2}\right) > 0 \tag{54}$$

(54) indicates that $t_{n_4 \to n_1} < 0$ which contradicts with $t_{n_4 \to n_1} > 0$. Therefore, such subtour does not exist, which completes the proof.